\documentclass[hidelinks,onefignum,onetabnum]{siamart251216}
\usepackage{amssymb}
\usepackage{amsfonts}
\usepackage{graphicx}
\usepackage[caption=false]{subfig}
\usepackage{url}
\makeatletter\g@addto@macro\UrlBreaks{\do\-}\makeatother
\usepackage{algpseudocode}
\usepackage{multirow}
\usepackage{colortbl}
\usepackage{array, booktabs}
\usepackage{longtable}
\makeatletter
\def\LT@makecaption#1#2#3{%
  \LT@mcol\LT@cols c{\hbox to\z@{\hss\parbox[t]\LTcapwidth{%
    \def\@captype{table}#1{\@makecaption{#2}{#3}}\endgraf}%
  \hss}}}
\makeatother
\usepackage{xspace}
\usepackage[shortlabels]{enumitem}
\usepackage{listings}
\usepackage{bm}

\newcommand{\secref}[1]{\cref{#1}}

\newcommand{\fref}[1]{\cref{#1}}
\newcommand{\tref}[1]{\cref{#1}}

\renewcommand{\eqref}[1]{Eq.~(\ref{#1})}
\newcommand{\eref}[1]{(\ref{#1})}
\newcommand{\pfrac}[2]{\frac{\partial #1}{\partial #2}}
\newcommand{\ppfrac}[2]{\frac{\partial^2 #1}{\partial #2^2}}

\newcommand{\mat}[1]{\textup{\textrm{\textbf{#1}}}}

\definecolor{light}{rgb}{0.8,0.8,0.8}
\definecolor{medium}{rgb}{0.6,0.6,0.6}
\definecolor{dark}{rgb}{0.4,0.4,0.4}
\definecolor{darkmed}{rgb}{0.3,0.3,0.3}
\definecolor{darkest}{rgb}{0.2,0.2,0.2}
\definecolor{Black}{rgb}{0,0,0}
\definecolor{White}{rgb}{1,1,1}
\definecolor{lightpurple}{rgb}{0.78823,0.709803,0.74509}
\definecolor{lightpurpletext}{rgb}{0.788235,0.5529411,0.658823}
\definecolor{skyblue}{rgb}{0.80392,0.866666,0.92941}
\definecolor{skybluetext}{rgb}{0.61568627,0.7647058,0.913725}
\definecolor{darkgreen}{rgb}{0.3137254,0.458823,0.18431}
\definecolor{foliagegreen}{rgb}{0.188,0.415,0.105}
\definecolor{steelbluegrey}{rgb}{0.1961,0.2353,0.2392}
\definecolor{highlightblue}{rgb}{0.4078,0.6431,0.85}
\definecolor{matlabblue}{rgb}{0,0.2705,0.85}
\definecolor{darkred}{rgb}{0.8,0.1725,0}
\definecolor{fireenginered}{rgb}{0.505,0.1411,0}
\definecolor{darkpurple}{rgb}{0.6431,0.3137,0.8509}
\definecolor{gaylordpurple}{rgb}{0.416,0.204,0.549}
\definecolor{deludedorange}{rgb}{0.7409,0.4392,0}
\definecolor{darksalmon}{rgb}{0.9137,0.411,0.706}

\headers{Parallel iterative methods for ILU triangular solves}{S. Dargaville and C. C. Pain}

\title{A comparison of parallel iterative methods for ILU triangular solves\thanks{Submitted to the editors \today. \funding{This work was supported by internal funding from the Department of Earth Science and Engineering, Imperial College London.}}}

\author{Steven Dargaville\thanks{Applied Modelling and Computation Group, Imperial College London, SW7 2AZ, UK \mbox{(\email{dargaville.steven@gmail.com})}.}
\and Christopher C. Pain\footnotemark[2]}

\crefname{algorithm}{Algorithm}{Algorithms}
\Crefname{algorithm}{Algorithm}{Algorithms}
\crefname{section}{Section}{Sections}
\crefname{subsection}{Section}{Sections}
\crefname{subsubsection}{Section}{Sections}
\Crefname{section}{Section}{Sections}
\Crefname{subsection}{Section}{Sections}
\Crefname{subsubsection}{Section}{Sections}

\ifpdf
\hypersetup{
  pdftitle={A comparison of parallel iterative methods for ILU triangular solves},
  pdfauthor={S. Dargaville and C. C. Pain}
}
\fi

\begin{document}
\maketitle

\begin{abstract}
Exact triangular solves parallelise poorly, limited by the size and number of the level sets, so on GPUs the factors from incomplete LU (ILU) factorisations are often solved approximately with parallel iterative methods. We compare four such methods: a Jacobi method, an Incomplete Sparse Approximate Inverse (ISAI), GMRES polynomials, and AIRG, a reduction multigrid built on GMRES polynomials, for which we introduce a parallel coarse/fine splitting that bounds the diagonal dominance ratio of the fine-fine block. We test them with ILU(0) and ILU(1) factors of canonical PDE problems and 91 SuiteSparse/Matrix Market matrices, the latter on a single GPU. AIRG is the most robust on the PDE factors, with nearly flat iteration counts under refinement and no per-problem tuning, but has the largest setup and memory costs. GMRES polynomials are more robust than the Jacobi method at a similar cost. As the inner solve of an ILU preconditioner with tuned inner iteration counts, ISAI gives the largest median speedups over a hipSPARSE exact triangular solve (4.5$\times$ and 6.8$\times$ for ILU(0) and ILU(1)), though this was dominated by smaller problems; on the largest ILU(1) problems all four methods give similar median speedups of 2.0--2.6$\times$.
\end{abstract}
\section{Introduction}
\label{sec:Introduction}
Incomplete LU (ILU) factorisations are often used as smoothers in multigrid methods or directly as preconditioners for difficult, sparse linear systems. ILU factorisations decompose a matrix $\mat{A}\approx \mat{L} \mat{U}$ into two triangular factors while limiting fill-in (the setup phase). Hence an approximate solution to the linear system $\mat{A}\mat{x}=\mat{b}$ can be found with $\mat{x} \approx \mat{U}^{-1}\mat{L}^{-1}\mat{b}$ (the solve phase). Parallel algorithms for computing both phases are vital.

The setup phase can be parallelised by level-scheduling \cite{Saad2003}, which searches for sets of rows/columns that can be calculated independently of each other, or by the nonlinear fixed-point iteration of \cite{Chow2015} (``ParILU''), which computes each entry of $\mat{L}$ and $\mat{U}$ independently and is now widely used on GPUs. The solve phase then requires the solution of two sparse triangular systems (SpTRSV). In serial this is very efficient, as forward/backward substitution (i.e., a Gauss--Seidel method in the triangular ordering) applies $\mat{U}^{-1}$ and $\mat{L}^{-1}$ exactly, but Gauss--Seidel methods are poorly suited to a parallel environment. Level-scheduling can be used \cite{Naumov2011, Saad2003}, as can synchronisation-free \cite{Liu2016a} or block/tiling algorithms \cite{Lu2020}. Level-scheduling can have limited parallelism, as any one level set may not have enough unknowns to saturate the available parallelism and unknowns in other level sets must wait for updated unknowns to reach them; in general, exact solves of problems with triangular structure are not scalable in parallel.

Given this and the fact that ILU methods only produce approximate factorisations, parallel iterative methods with large tolerances are often used in the solve phase. $\mat{L}$ and $\mat{U}$ are nonsymmetric and often strongly non-normal; there is limited convergence theory available for iterative methods applied to such systems. It has been observed, however, that the triangular factors are often close to diagonally dominant for stable incomplete factorisations. Given this, a few iterations of a Jacobi method or sparse approximate inverses (SAIs) have proven very effective, as their application is very parallel and suited to GPUs \cite{Chow2015, Anzt2015, Chow2018, Thomas2021, Anzt2018a, Lu2021, Thomas2023}, giving large speedups over exact triangular solves. Vendor sparse triangular solves continue to improve, however; on an AMD MI250X, \cite{Thomas2023} found the gain from iterative triangular solves in the ILU smoother of an AMG method was around 1.5$\times$.

\begin{sloppypar}
While parallel implementations of algebraic multigrids are readily available, they have seen limited use replacing exact triangular solves as their convergence often degrades in nonsymmetric systems. In the last decade, however, reduction multigrids using AIR (Approximate Ideal Restriction) have been developed with excellent performance in nonsymmetric systems \cite{Southworth2017, Manteuffel2019, Manteuffel2019a, Ali2024, Zaman2024, Dargaville2024a, Dargaville2024}. Recently we showed excellent weak scaling up to 11k GPUs when solving triangular (and near-triangular) systems from discretisations of advection problems with AIRG (Approximate Ideal Restriction with GMRES polynomials) \cite{Dargaville2026}. 

In this work we compare four different parallel iterative methods for solving the triangular factors from ILU factorisations: a Jacobi method, an Incomplete Sparse Approximate Inverse (ISAI), GMRES polynomials, and AIRG. We test their robustness on difficult problems and their performance on a GPU, both when solving the factors alone and when used within an ILU preconditioner. We restrict ourselves to incomplete factorisations, as the factors from complete LU factorisations are typically much denser.

The main contributions of this paper are: (i) a new CF splitting algorithm for AIRG on GPUs, which produces a fine-fine block with a user-specified maximum diagonal dominance ratio (\secref{sec:cf_splitting}); (ii) the use of AIRG and GMRES polynomials to solve ILU triangular factors; and (iii) a comparison of the convergence, parallelism and GPU performance of these methods (\secref{sec:Results}). The new CF splitting and all of the iterative methods compared, along with cross-platform GPU implementations, are available in the open-source \textsc{PFLARE}\footnote{\url{https://github.com/PFLAREProject/PFLARE}} library, which can be installed directly through the PETSc configure \cite{petsc-web-page} with \texttt{--download-pflare}.
\end{sloppypar}
\section{ILU factorisations}
\label{sec:ilu}
The quality of an ILU factorisation as a preconditioner depends heavily on the ordering of the unknowns \cite{Benzi1999a}, the fill-in allowed and any pivoting used; Chow and Saad \cite{Chow1997} give an excellent characterisation of the properties that contribute to that quality. These choices affect the condition number, row scaling, diagonal dominance and normality of the resulting factors, and a poor ordering can produce $\mat{L}$ and $\mat{U}$ factors that are essentially unsolvable, sometimes even by direct methods. Our aim is not to examine how to produce good ILU factorisations, but rather how to iteratively solve the resulting factors (for nonsymmetric systems, how to reliably produce ``good'' approximate factorisations remains an open question), so we take several steps to improve the quality of $\mat{L}$ and $\mat{U}$.

Before factorising, we left-scale the matrix by the inverse diagonal. We then loop through a number of different orderings and, for each, produce an ILU(0) or ILU(1) factorisation and use it as a preconditioner for GMRES(30) applied to the original system $\mat{A}\mat{x}=\mat{b}$, with an exact triangular solve used to apply the factors. We select the ordering with the lowest iteration count among those that converge within 2000 iterations or, if none converge, the one with the lowest final relative residual, and use it for all subsequent solves; we note the winning ordering in the results. The orderings we try are: \textsc{natural}, \textsc{nd} (nested dissection), \textsc{rcm} (reverse Cuthill--McKee), \textsc{qmd} (quotient minimum degree), \textsc{1wd} (one-way dissection) and \textsc{rowlength} (orderings from PETSc), and \textsc{colamd} (column approximate minimum degree), \textsc{mmd\_ata} and \textsc{mmd\_at\_plus\_a} (multiple minimum degree on $\mat{A}^\textrm{T}\mat{A}$ and $\mat{A}^\textrm{T}+\mat{A}$), which are fill-reducing orderings from SciPy (see the start of \cref{sec:sm_results} of the supplementary materials for how these are extracted). Every ordering is applied symmetrically.

AIRG requires access to assembled sparse matrices, so we compute $\mat{L}$ and $\mat{U}$ explicitly with an exact level-based incomplete factorisation, ILU($k$); on the GPU we use the Kokkos Kernels \cite{Rajamanickam2021, Trott2022} \texttt{spiluk} implementation, which parallelises the factorisation with level-scheduling. The factorisation returns $\mat{L}$ as a unit lower-triangular matrix and $\mat{U}$ with the pivots on its diagonal, so we left-scale $\mat{U}$ by its inverse diagonal to give a unit diagonal, which typically improves the behaviour of the iterative solves applied to it (e.g., see \cite{Thomas2021} for an examination of different scaling approaches). We do not perform any pivoting, either when computing the factors or the orderings, so all of the matrices we test have nonzero diagonals, which the symmetric reorderings preserve; more robust ILU factorisations using pivoting are available. We now examine the iterative methods used to solve $\mat{L}$ and $\mat{U}$, written for a generic system $\mat{A}\mat{x}=\mat{b}$ unless a specific property of the triangular factors is used.
\section{ISAI}
\label{sec:isai}
Sparse Approximate Inverses (SAIs) \cite{Kolotilina1993, Grote1997, Chow1998, Chow2000} can be used to compute an approximation $\mat{M} \approx \mat{A}^{-1}$. A left inverse is computed by minimising $\|\mat{M} \mat{A} - \mat{I}\|_\textrm{F}^2$ over all $\mat{M}$ with a given sparsity pattern. This minimisation decouples into an independent least-squares problem for each row $j$ of $\mat{M}$ (or each column for a right inverse). We write $\mathcal{J}$ as the indices of the nonzero columns in the $j$th row of $\mat{A}$ and $\mathcal{I}$ as the ``shadow'' of $\mathcal{J}$ (i.e., the indices of the nonzero columns in the rows of $\mathcal{J}$) in $\mat{A}$; the least-squares problem for row $j$ then only involves the small block $\mat{A}(\mathcal{J}, \mathcal{I})$ and determines the nonzero entries $\mat{M}(j, \mathcal{J})$.

In this work we only use the ``incomplete'' SAI (ISAI) \cite{Anzt2018a}; the incomplete version replaces the shadow $\mathcal{I}$ with $\mathcal{J}$, so that the least-squares problem for each row becomes square and is solved exactly by
\begin{equation}
\mat{M}(j, \mathcal{J}) \mat{A}(\mathcal{J}, \mathcal{J}) = \mat{I}(j, \mathcal{J})
\label{eq:isai_small}
\end{equation}
for the entries of $\mat{M}$ in each row, $j$. An ISAI is typically cheaper to form than a SAI as $\mathcal{J} \subseteq \mathcal{I}$, and ISAIs can produce better preconditioners than SAIs \cite{Anzt2018a}.
We choose a sparsity pattern for $\mat{M}$ in advance, as this fixes the index sets $\mathcal{J}$ (and hence $\mathcal{I}$) used in each row solve. In this work we build an assembled $\mat{M}$ with the sparsity pattern of $\mat{A}$ itself, rather than that of a higher power of $\mat{A}$.
\section{Jacobi method/Neumann polynomial}
\label{sec:jacobi}
A simple iterative method for solving triangular factors is a Jacobi iteration. As noted in \secref{sec:ilu}, both $\mat{L}$ and $\mat{U}$ have been scaled to have a unit diagonal, so for a generic triangular factor $\mat{T} \in \{\mat{L}, \mat{U}\}$ we can write $\mat{T} = \mat{I} - \mat{G}$, where $\mat{G} = \mat{I} - \mat{T}$ is strictly triangular (strictly lower for $\mat{L}$, strictly upper for $\mat{U}$). Hence the Jacobi iteration for $\mat{T}\mat{x}=\mat{b}$ reduces to $\mat{x}^{(m+1)} = \mat{G}\,\mat{x}^{(m)} + \mat{b}$, with iteration matrix $\mat{G}$. Starting from a zero initial guess $\mat{x}^{(0)} = \mat{0}$, after $m$ iterations this gives
\begin{equation}
\mat{x}^{(m)} = \sum_{j=0}^{m-1} \mat{G}^j \, \mat{b} = p_{m-1}(\mat{T})\,\mat{b}, \qquad p_{m-1}(\mat{T}) = \mat{I} + \mat{G} + \ldots + \mat{G}^{m-1} \approx \mat{T}^{-1},
\label{eq:neumann}
\end{equation}
which is the application of a Neumann polynomial of degree $m-1$ \cite{Saad2003, Chow2015}.
We use the two descriptions interchangeably in what follows. The Jacobi iteration requires only sparse matrix-vector products (SpMVs), so it is matrix-free, highly parallel and well suited to GPUs; this is the primary reason a few iterations of a Jacobi method have become a popular approximate triangular solve for ILU factorisations \cite{Chow2015, Anzt2015, Chow2018, Thomas2021}. The truncated series can also be assembled explicitly, with dropping, as a type of sparse approximate inverse \cite{Lu2021}.

Because $\mat{G}$ is strictly triangular it is nilpotent, with $\mat{G}^n = \mat{0}$ for a factor of size $n$, so the Neumann series in \eref{eq:neumann} terminates exactly after $n$ terms and the spectral radius of the iteration matrix is zero. Asymptotic convergence is therefore trivially guaranteed, with the iteration exact after at most $n$ steps. This is misleading in practice, however, as we only ever apply a low-degree polynomial with $m \ll n$ and a zero spectral radius says nothing about the size of $\mat{G}^m$ for small $m$. When $\mat{G}$ is non-normal, as is typically the case for the factors of an ILU factorisation, $\| \mat{G}^m \|$ can grow by many orders of magnitude before eventually decaying to zero, so a low-order Neumann polynomial can be highly inaccurate and the iteration can appear to diverge. The practical convergence is governed by the degree of non-normality of $\mat{G}$ and is better characterised by its pseudospectra or field of values than by its spectrum \cite{Trefethen2005}.
\section{GMRES polynomial}
\label{sec:gmres}
Rather than fixing the coefficients to those of the Neumann series in \eref{eq:neumann} (which are all unity), a GMRES polynomial chooses them to minimise a residual. The GMRES polynomial of degree $m-1$ written in a monomial basis (shown here for simplicity; we do not use this form as it is not numerically stable at high order) is
\begin{equation}
q_{m-1}(\mat{A}) = \alpha_0 \mat{I} + \alpha_1 \mat{A} + \alpha_2 \mat{A}^2 + \ldots + \alpha_{m-1} \mat{A}^{m-1} \approx \mat{A}^{-1},
\label{eq:gmres_poly}
\end{equation}
where the coefficients $\bm{\alpha}$ are those implicitly constructed by $m$ steps of GMRES applied to $\mat{A}\mat{x} = \mat{b}$; that is, $\mat{x} \approx q_{m-1}(\mat{A})\mat{b}$ is the GMRES iterate that minimises the residual $\| \mat{b} - \mat{A}\, q_{m-1}(\mat{A}) \mat{b} \|_2$ over the Krylov subspace \cite{Nachtigal1992, Loe2022, Embree0}. We compute the GMRES polynomial in a Newton basis \cite{Loe2022} once during the setup phase, from a single random right-hand side, and then reuse the same fixed polynomial for every subsequent application. This produces a stationary polynomial that can be applied matrix-free with only SpMVs (and no dot products), so, like the Neumann polynomial, it maps well onto GPUs.
\section{AIRG}
\label{sec:reduction}
AIRG is a reduction multigrid built using the GMRES polynomials described in \secref{sec:gmres}; please see \cite{Southworth2017, Manteuffel2019, Manteuffel2019a, Zaman2024, Dargaville2024a} for further details. An important feature of AIRG, in the context of the present work, is that it does not require a specific ordering of the unknowns, or any block or lower-triangular structure in the matrix, to be applied. Given this, all the components of AIRG are well suited to GPUs. The starting point of a reduction multigrid is a splitting of the unknowns into a set of ``coarse'' and ``fine'' points, with $n_\textrm{C} + n_\textrm{F} = n$. This coarse/fine (CF) splitting allows us to write the block system
\begin{equation}
\begin{bmatrix}
\mat{A}_\textrm{ff} & \mat{A}_\textrm{fc} \\
\mat{A}_\textrm{cf} & \mat{A}_\textrm{cc}
\end{bmatrix}
\begin{bmatrix}
\mat{x}_\textrm{f} \\
\mat{x}_\textrm{c}
\end{bmatrix} =
\begin{bmatrix}
\mat{b}_\textrm{f} \\
\mat{b}_\textrm{c}
\end{bmatrix},
\label{eq:air_two}
\end{equation}
and form prolongation and restriction operators of the form
\begin{equation}
\mat{P} =
\begin{bmatrix}
\mat{W} \\
\mat{I}
\end{bmatrix}, \quad
\mat{R} =
\begin{bmatrix}
\mat{Z} & \mat{I}
\end{bmatrix}.
\label{eq:prolong}
\end{equation}
The ideal restrictor $\mat{Z}=-\mat{A}_\textrm{cf} \mat{A}_\textrm{ff}^{-1}$ enforces $\mat{R}\mat{A}=\mat{0}$ on F-points, with the ideal prolongator $\mat{W}=-\mat{A}_\textrm{ff}^{-1} \mat{A}_\textrm{fc}$ enforcing $\mat{A}\mat{P}=\mat{0}$ on F-points. The action of $\mat{A}_\textrm{ff}^{-1}$ is in general dense and cannot be formed exactly. A practical reduction multigrid therefore replaces it with a sparse approximation $\hat{\mat{A}}_\textrm{ff}^{-1} \approx \mat{A}_\textrm{ff}^{-1}$ and forms approximate ideal grid-transfer operators
\begin{equation}
\mat{Z} = -\mat{A}_\textrm{cf} \hat{\mat{A}}_\textrm{ff}^{-1}, \quad \mat{W} = -\hat{\mat{A}}_\textrm{ff}^{-1} \mat{A}_\textrm{fc}.
\label{eq:approx_ideal}
\end{equation}
AIR methods typically only use $\mat{Z}$ in \eref{eq:approx_ideal} and form a simpler approximation for $\mat{W}$ \cite{Manteuffel2019}; we use a classical one-point prolongator, in which each F-point is interpolated with unit weight from its strongest C-point neighbour (the largest-magnitude entry in its row of $\mat{A}_\textrm{fc}$). The coarse grid matrix is computed as $\mat{A}_\textrm{coarse}=\mat{R}\mat{A}\mat{P}$ (which is an approximation to the Schur complement $\mat{A}_\textrm{cc}-\mat{A}_\textrm{cf}\mat{A}_\textrm{ff}^{-1}\mat{A}_\textrm{fc}$) and the procedure is applied recursively to build the multigrid hierarchy. To limit fill, entries of $\mat{Z}$ and of each coarse grid matrix smaller than a relative drop tolerance ($10^{-4}$ and $10^{-6}$, respectively) times the largest-magnitude entry in their row are removed.

AIRG uses the GMRES polynomial of \secref{sec:gmres} applied to the fine-fine block as the approximate inverse, $\hat{\mat{A}}_\textrm{ff}^{-1} \approx q_{m-1}(\mat{A}_\textrm{ff}) \approx \mat{A}_\textrm{ff}^{-1}$. This ensures that using AIRG in the solve phase only requires SpMVs. Note that approximations to $\mat{A}_\textrm{ff}^{-1}$ play two distinct roles in the multigrid: they are used to build the approximate ideal restrictor in \eref{eq:approx_ideal}, and also as an F-point smoother on the up cycle, with no C-point or down-smoothing required \cite{Dargaville2024a}. As in \cite{Dargaville2024a, Dargaville2025}, we use a single assembled matrix representation of $q_{m-1}(\mat{A}_\textrm{ff})$ for both, with fill-in controlled by enforcing the sparsity of $\mat{A}_\textrm{ff}$ on each of the matrix powers; the only polynomial we apply matrix-free is the higher-order coarse-grid solver described in \secref{sec:independent}.

Good performance of AIRG with low-order GMRES polynomials depends on the CF splitting producing an $\mat{A}_\textrm{ff}$ block that is amenable to a good sparse approximate inverse, which in practice means selecting fine rows that are as diagonally dominant as possible. We first place this requirement in the context of other multilevel methods, before detailing the CF splitting algorithm we use to ensure it.
\subsection{Relationship to multilevel ILU}
\label{sec:multilevel_ilu}
AIRG and other reduction multigrids are related to, but distinct from, multilevel ILU methods \cite{Saad1996, Saad2002, Saad2005, Vassilevski2008}. Like reduction multigrids, multilevel ILU methods build a hierarchy by eliminating a leading submatrix on each level and deferring the remaining coupling to a Schur complement that is treated recursively. A multilevel ILU is an approximate block LDU factorisation and so applies its inner solve additively, within the block-diagonal $\mat{D}$ factor, whereas a reduction multigrid applies its F-point smoothing multiplicatively, after the coarse-grid correction. Obtaining parallelism in any of these methods typically begins with an independent set, and the canonical example is level-scheduling. For an ILU factorisation, level-scheduling is applied directly to the triangular factors $\mat{L}$ and $\mat{U}$, parallelising both their setup and solve phases. Consider a lower-triangular factor $\mat{L}$: if we symmetrically permute $\mat{L}$ (i.e., a reordering) so that each level set is contiguous, $\mat{L}$ becomes block lower-triangular and, crucially, each diagonal block is itself diagonal, because the unknowns within a single level set do not reference one another. The unknowns in each block are those of one level set; the blocks must be processed sequentially (given the block lower-triangular structure), but all the unknowns within a block can be solved in parallel.

Viewed as a multilevel method, level-scheduling is one particular way of producing a diagonal $\mat{A}_\textrm{ff}$ on every level. Any independent set of unknowns gives a diagonal $\mat{A}_\textrm{ff}$, as its members share no couplings; in general there are many such sets. The set used by level-scheduling is the maximal choice for which $\mat{A}_\textrm{fc}=\mat{0}$, with its fine rows being the sources of the dependency graph and so having no sub-diagonal entries. The Schur complement $\mat{A}_\textrm{cc}-\mat{A}_\textrm{cf}\mat{A}_\textrm{ff}^{-1}\mat{A}_\textrm{fc}$ then reduces to $\mat{A}_\textrm{cc}$, so no fill is produced and the reduced matrix remains triangular. Recursing reproduces level-scheduling exactly, i.e., the block Gaussian elimination of the forward substitution sweep.

Level-scheduling sits at one end of a family of multilevel methods built on a diagonal $\mat{A}_\textrm{ff}$, which trade fill against the number of levels. Enforcing $\mat{A}_\textrm{fc}=\mat{0}$ restricts each level to the current sources of the graph and so leaves as many levels as the critical path. Relaxing that requirement, keeping only the diagonal $\mat{A}_\textrm{ff}$ of an independent set whose members need not be sources, permits larger F-blocks and fewer levels at the cost of fill in the Schur complement. Fewer levels than level-scheduling can also be had without fill in some cases; two examples are the partitioned inverse of \cite{Alvarado1993} and cyclic reduction \cite{Golub1992}. Alvarado and Schreiber \cite{Alvarado1993} group consecutive columns of $\mat{L}$ into factors that are each invertible without fill, so that $\mat{L}^{-1}$ is applied exactly as a product of a few sparse factors. Because it is exact, how far the level sets can be merged is fixed by the sparsity of $\mat{L}$, and the merged product form is only conditionally stable \cite{Higham1994}. Cyclic reduction \cite{Golub1992} uses a maximal independent set on each level, and remains fill-free when applied to tridiagonal systems. In general, larger $\mat{A}_\textrm{ff}$ blocks and hence fewer levels allow more parallelism.

ILUM (ILU with multi-elimination) and its block version BILUM \cite{Saad1996, Saad1999a} use a maximal independent set on matrices of general sparsity: the independent set still gives a diagonal/block-diagonal $\mat{A}_\textrm{ff}$, which is inverted trivially and exactly, but eliminating it now produces fill in the Schur complement that must be controlled by dropping. Later multilevel ILU methods such as \cite{Saad2005} abandon diagonal $\mat{A}_\textrm{ff}$ in search of even larger F-blocks and fewer levels. They instead select a non-diagonal $\mat{A}_\textrm{ff}$ whose inverse is well behaved, in practice a diagonally dominant block. Such an $\mat{A}_\textrm{ff}$ is no longer trivial to invert exactly, and so its inverse is itself only approximated, by an (incomplete) factorisation of $\mat{A}_\textrm{ff}$.

Both of these steps move the hierarchy away from an exact solve: dropping the fill makes the elimination incomplete and approximating $\mat{A}_\textrm{ff}^{-1}$ makes the inversion inexact. The relevant distinction for parallelism is therefore not whether the solve on each level is exact, but how the (inexact) inverse of $\mat{A}_\textrm{ff}$ is applied. In these multilevel ILU methods the approximate $\mat{A}_\textrm{ff}^{-1}$ is built from (incomplete) triangular factors and applied by a triangular solve, so the parallelism available on each level is once again governed by the sparsity structure and ordering of those factors, reintroducing the level-set problem we began with. We also note that these multilevel ILU methods are built for and applied to the original matrix $\mat{A}$; each is a self-contained multilevel preconditioner for $\mat{A}$, rather than a means of solving an already-computed triangular factor. Multilevel ILU methods could be used to solve already-computed triangular factors, but this would mean computing a further ILU factorisation of $\mat{L}$ or $\mat{U}$ and a triangular solve on every level, which is both redundant and exposes less parallelism.

In contrast, AIRG applied to either the original matrix $\mat{A}$ or to the triangular factors $\mat{L}$ and $\mat{U}$ typically retains more parallelism. It does this by coarsening based on the same idea as in the multilevel ILU of \cite{Saad2005}, targeting the extraction of a diagonally dominant non-diagonal $\mat{A}_\textrm{ff}$ (with an algorithm detailed in \secref{sec:cf_splitting}), but then applying the approximate inverse of $\mat{A}_\textrm{ff}$ as a GMRES polynomial rather than by a triangular solve. The essential point of AIRG (or of a comparable multilevel LDU/Schur-complement scheme built on a GMRES polynomial approximation of $\mat{A}_\textrm{ff}^{-1}$) is that at no point does the parallelism on any level depend on the ordering of the unknowns. The price of this is robustness; multilevel ILU factorisations applied to the original matrix $\mat{A}$ are typically very robust, as they can utilise much of the literature around producing good ILU factorisations in difficult problems, whereas convergence theory and practical guidance for GMRES (and therefore GMRES polynomials) is less developed. In this work we focus on the SpTRSV problem: ILU factorisations are widely used, robust and readily available on GPUs, so we take the $\mat{L}$ and $\mat{U}$ factors they have already (robustly) produced as given and use AIRG to solve them. We leave a study of applying AIRG to $\mat{A}$ directly in these problems for future work. We now consider a parallel CF splitting algorithm for producing diagonally dominant $\mat{A}_\textrm{ff}$.
\subsection{DIAG\_DOM CF splitting}
\label{sec:cf_splitting}
In our previous work \cite{Dargaville2025} we introduced the PMISR DDC algorithm (Parallel Modified Independent Set for Reduction multigrids and Diagonal Dominance Cleanup) in order to produce diagonally dominant $\mat{A}_\textrm{ff}$ blocks. The first pass, PMISR, is a slight modification of PMIS (the Parallel Modified Independent Set algorithm of \cite{Sterck2006}) that computes an independent set in $\mat{S}+\mat{S}^\textrm{T}$ (where $\mat{S}$ is the classical strength-of-connection matrix, in which an off-diagonal entry $A_{ij}$ is strong if $|A_{ij}| \ge \theta \max_{k \neq i} |A_{ik}|$ for a strength threshold $\theta$) such that no large off-diagonal entries appear in the resulting $\mat{A}_\textrm{ff}$. The second pass, DDC, then flips the F rows with the worst diagonal dominance back to C. We define the diagonal dominance ratio of an F row $i$ restricted to the F columns as
\begin{equation}
r_\textrm{dd}(i) = \frac{\sum_{j \in \mathcal{F},\, j \neq i} |A_{ij}|}{|A_{ii}|}, \qquad \mathcal{F}=\{i : c(i) = \texttt{F}\},
\label{eq:dd_ratio}
\end{equation}
where $\mat{c}$ is the vector of C/F markers. The DDC pass takes a user-specified local fraction $f$ and swaps at least the $f n_\textrm{F}$ worst F rows (by $r_\textrm{dd}$) back to C; in our implementation it also always swaps every F row with $r_\textrm{dd} \ge 1$, so the resulting $\mat{A}_\textrm{ff}$ is always diagonally dominant.

For discretisations of advection-dominated PDEs we have found this gives an $\mat{A}_\textrm{ff}$ that is sufficiently diagonally dominant for AIRG to converge well. PMISR DDC, however, gives no control over how diagonally dominant $\mat{A}_\textrm{ff}$ is beyond the fixed bound of one, while the fraction $f$ of additional rows flipped is a problem-dependent parameter. This is suboptimal for the triangular factors considered in this work, where a stronger bound can be needed (\secref{sec:comparison}). To improve robustness we therefore introduce a variant of PMISR DDC that iterates to a user-specified target $\tau$ on the maximum diagonal dominance ratio across all rows of $\mat{A}_\textrm{ff}$; we refer to this as the DIAG\_DOM splitting. Selecting the largest subset of rows of a matrix whose diagonal dominance ratios are below a given threshold is NP-complete \cite{MacLachlan2007} and hence we do not pursue an exact solution. Heuristics for determining diagonally dominant submatrices have been investigated within the context of different multilevel methods \cite{Saad2002, Saad2005, MacLachlan2007, MacLachlan2007a, Zaman2023, Dargaville2025, Dargaville2026}.

We follow on from this previous work directly and seek a parallel algorithm that avoids any inherently sequential selection. As such we only use maximal independent set algorithms as our building blocks. An outline of our algorithm is to start with a first pass which selects candidate F-points and then rely on a number of iterations after this to reject ``bad'' F-points. The two phases pursue different objectives, and our measure choices reflect this.

The first pass selects on connectivity in $\mat{S}+\mat{S}^\textrm{T}$ rather than on $r_\textrm{dd}$ directly. Typically we find the hard loss of diagonal dominance is driven by clusters of mutually large, comparable off-diagonal entries competing for the same diagonal budget; the connectivity graph of $\mat{S}+\mat{S}^\textrm{T}$ encodes exactly this competition, so a maximal independent set thins each cluster in a single pass while also favouring the low-connectivity rows whose elimination keeps the coarse-grid matrix sparse. A measure based on $r_\textrm{dd}$ achieves neither, and we found that using $r_\textrm{dd}$ as the first-pass measure substantially increased the number of rejection iterations.

Hence our DIAG\_DOM CF splitting consists of a first pass of PMISR with $\mat{S}+\mat{S}^\textrm{T}$, where $\mat{S}$ is formed from $\mat{A}$ with the strength threshold set to $\theta=\tau$, the target diagonal dominance ratio, so that DIAG\_DOM has a single parameter. The rejection step then iterates as many times as required, converting F-points to C until the resulting $\mat{A}_\textrm{ff}$ meets the target $\tau$. Swapping rows one at a time would trivially allow this, but would not be sufficiently parallel.

Instead, on each iteration we compute an independent set over the F rows that do not yet satisfy $r_\textrm{dd} < \tau$, using the same PMISR routine as in the first pass on $\mat{A}_\textrm{ff}+\mat{A}_\textrm{ff}^\textrm{T}$, and swap its members to C. This swaps many non-adjacent F-points simultaneously, but ensures that for any pair of strongly coupled F-F points, only one is swapped per iteration; the neighbours whose $r_\textrm{dd}$ improves as a result are re-evaluated on the next iteration. As PMISR selects the rows with the \textit{smallest} measure, we invert $r_\textrm{dd}$ and use the measure
\begin{equation}
\mu(i) \;=\; M \,-\, r_\textrm{dd}(i) \,+\, \delta\, \xi(i),
\label{eq:measure}
\end{equation}
where $\xi(i) \in [0,1)$ is a random value used for tie-breaking and $\delta=10^{-10}$ is small enough not to change the relative ordering of rows with distinct $r_\textrm{dd}$. The constant $M \ge \max_i r_\textrm{dd}(i) + 1$ keeps $\mu(i) \ge 1$, as PMISR would otherwise classify a row with $|\mu| < 1$ (normally a row with no strong neighbours) directly as F. Rows already satisfying $r_\textrm{dd} < \tau$ are assigned $\mu=+\infty$ and are never swapped. We also tried measuring strength relative to the diagonal in the first pass, $|A_{ij}| \ge \tau |A_{ii}|$, and scaling $r_\textrm{dd}$ by the inverse of the number of F-F neighbours to keep $\mat{A}_\textrm{ff}$ sparser \cite{Saad2005}; the former gave a worse first pass that required considerably more work in the rejection step, and the latter gave no large improvement, as our first pass already targets this objective.

The freeze of already-satisfied rows ensures the rejection iterations terminate. Flipping any F row to C only ever removes terms from another row's F-F sum in \eref{eq:dd_ratio}, so $r_\textrm{dd}$ is monotone non-increasing for every surviving row; a row that already satisfies the target therefore stays satisfied. Each iteration strictly decreases the number of rows violating the target, so the loop cannot continue indefinitely and terminates once every row meets $\tau$. In the worst case this termination may take as many iterations as there are F-points produced by the first pass, but in practice we find that for sparse factors only a few iterations are required.

This is summarised in \cref{alg:pmisr_dd}. The use of an implicit transpose lets us pass $\mat{A}_\textrm{ff}$ directly to the independent set algorithm without ever forming $\mat{A}_\textrm{ff}+\mat{A}_\textrm{ff}^\textrm{T}$ explicitly. The whole algorithm is well suited to GPUs: the diagonal dominance ratio is embarrassingly parallel to compute and the only global reduction required is a single max over $r_\textrm{dd}$ across all F rows per iteration. Each iteration involves computing an independent set, for which efficient GPU implementations are available \cite{Naumov2015, Falgout2021}.
\begin{algorithm}
\caption{DIAG\_DOM CF splitting iteration.}
\label{alg:pmisr_dd}
\begin{algorithmic}[1]
\Require Matrix $\mat{A}$, initial C/F markers $\mat{c}$ from first-pass PMISR on $\mat{S}+\mat{S}^{\textrm{T}}$, target ratio $\tau$
\Loop
   \State Extract $\mat{A}_\textrm{ff}$ from $\mat{A}$ using $\mat{c}$ \Comment{F rows and F columns of $\mat{A}$}
   \State Compute $r_\textrm{dd}(i)$ for every $i \in \mathcal{F}$ \Comment{\eqref{eq:dd_ratio}}
   \State $r_\textrm{max} \gets \max_{i \in \mathcal{F}} r_\textrm{dd}(i)$ over all ranks \Comment{Global reduction}
   \State \textbf{if} $r_\textrm{max} < \tau$ \textbf{then exit} \Comment{$\mat{A}_\textrm{ff}$ meets the target}
   \State Initialise swap markers $\mat{s}(i) \gets \texttt{candidate}$ for every $i \in \mathcal{F}$
   \For{$i \in \mathcal{F}$ \textbf{in parallel}}
      \If{$r_\textrm{dd}(i) < \tau$}
         \State $\mu(i) \gets +\infty$; \, $\mat{s}(i) \gets \texttt{keep}$ \Comment{Already satisfied, never swapped}
      \Else
         \State $\mu(i) \gets M - r_\textrm{dd}(i) + \delta \, \xi(i)$ \Comment{\eqref{eq:measure}}
      \EndIf
   \EndFor
   \State Run PMISR on $\mat{A}_\textrm{ff}+\mat{A}_\textrm{ff}^\textrm{T}$ with measure $\bm{\mu}$ \Comment{Largest $r_\textrm{dd}$ first}
   \State $\mat{s}(i) \gets \texttt{swap}$ for every candidate $i$ in the resulting independent set
   \For{$i \in \mathcal{F}$ \textbf{in parallel}}
      \If{$\mat{s}(i) = \texttt{swap}$}
         \State Flip $\mat{c}(i)$ from F to C
      \EndIf
   \EndFor
\EndLoop
\end{algorithmic}
\end{algorithm}
\subsection{Comparison of CF splittings}
\label{sec:comparison}
In this section we show that the DIAG\_DOM CF splitting in \cref{alg:pmisr_dd} can enforce a maximum diagonal dominance ratio while retaining a similar number of F-points to PMISR DDC. For many problems derived from PDEs, PMISR DDC can produce a sufficiently diagonally dominant block cheaply (as there is only one ``rejection'' step) and hence AIRG converges well. For more difficult problems, however, we find that DIAG\_DOM is more robust. For DIAG\_DOM we set $\tau=0.5$. For the PMISR DDC algorithm we use the same first pass as DIAG\_DOM, but for the second pass we do a single DDC iteration with $f=0.1$; this converts at least the worst 10\% of F-points to C-points, measured by $r_\textrm{dd}$, along with every F-point with $r_\textrm{dd} \ge 1$.

We show this for an example test matrix \textsc{bcsstk24} from the SuiteSparse collection. First we reorder the matrix with an \textsc{rcm} ordering and produce an ILU(0) factorisation of the matrix. We then run the CF splitting algorithms on the $\mat{L}$ factor (similar results are found for $\mat{U}$); this gives the $\mat{A}_\textrm{ff}$ that would be used on the top level of the AIRG multigrid applied to $\mat{L}$. For consistency we refer to the submatrix as $\mat{A}_\textrm{ff}$ rather than $\mat{L}_\textrm{ff}$.

\begin{figure}[thp]
\centering
\includegraphics[width=0.7\textwidth]{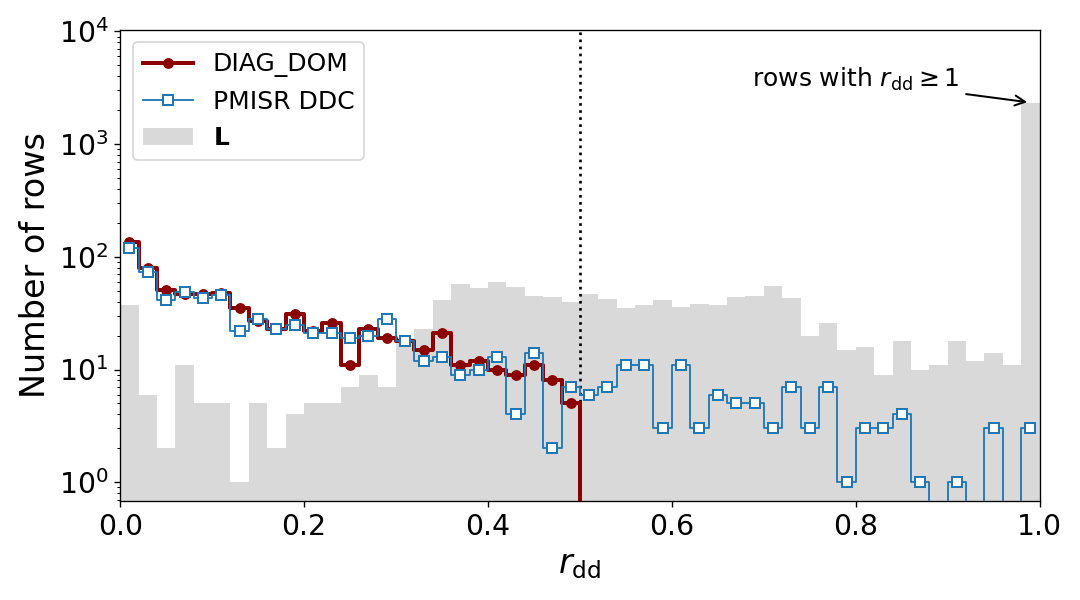}
\caption{Binned diagonal dominance ratio, $r_\textrm{dd}$, for each row in the lower-triangular factor $\mat{L}$ produced from an ILU(0) factorisation of the matrix \textup{\textsc{bcsstk24}}. The maximum $r_\textrm{dd}$ across any row in $\mat{L}$ is 11137 (i.e., $\mat{L}$ is far from diagonally dominant), but all rows with $r_\textrm{dd} \ge 1$ have been added to the top bin in the figure. The ratio is also shown for the $\mat{A}_\textrm{ff}$ submatrices resulting from a PMISR DDC and DIAG\_DOM CF splitting on $\mat{L}$.}
\label{fig:cf_compare}
\end{figure}

\fref{fig:cf_compare} shows that the matrix $\mat{L}$ is not diagonally dominant, with many rows having $r_\textrm{dd}>1$. For this matrix, the PMISR DDC algorithm results in a diagonally dominant $\mat{A}_\textrm{ff}$ matrix, with a maximum $r_\textrm{dd}$ equal to 0.998. The DIAG\_DOM algorithm, however, gives an $\mat{A}_\textrm{ff}$ matrix which is much more diagonally dominant; the maximum $r_\textrm{dd}$ is just below 0.5 (0.500 to three decimal places), which matches the maximum provided threshold $\tau=0.5$. Both CF splittings produce $\mat{A}_\textrm{ff}$ matrices of similar sizes for this matrix. PMISR DDC gives 788 F-points, with DIAG\_DOM giving 744 F-points; this is approximately 21--22\% of the unknowns (of $n=3562$). The difference shows in the convergence of AIRG, which applies the CF splitting on every level of the hierarchy. On the ILU(0) factors of \textsc{bcsstk24} AIRG takes 5 iterations on both $\mat{L}$ and $\mat{U}$ with either splitting, but on the harder $\mat{L}$ factor from ILU(1) it takes 18 iterations with PMISR DDC, compared with 6 with DIAG\_DOM (both take 5 on $\mat{U}$).

\tref{tab:tau_sweep} shows the effect of sweeping $\tau$ for the DIAG\_DOM splitting applied to the same $\mat{L}$ factor of \textsc{bcsstk24}. In every case the achieved maximum $r_\textrm{dd}$ matches the requested target $\tau$ to within a tight tolerance. As $\tau$ is decreased the resulting $\mat{A}_\textrm{ff}$ becomes smaller, ranging from around 28\% of the unknowns at $\tau=0.9$ down to around 12\% at $\tau=0.1$. A larger $\tau$ therefore produces a larger $\mat{A}_\textrm{ff}$, but the resulting block is less diagonally dominant and hence harder to approximate, which degrades the convergence of the resulting AIRG. Conversely, a smaller $\tau$ gives a more diagonally dominant (and more easily inverted) $\mat{A}_\textrm{ff}$, but extracts fewer F-points on each level. We discuss the implications of this trade-off for parallelism below.
\begin{table}[ht]
\centering
\caption{Effect of varying the target diagonal dominance ratio $\tau$ in the DIAG\_DOM CF splitting, applied to the $\mat{L}$ factor from an ILU(0) factorisation of \textup{\textsc{bcsstk24}}. The achieved maximum $r_\textrm{dd}$ is rounded to three decimal places; every value is strictly below its $\tau$.}
\label{tab:tau_sweep}
\begin{tabular}{cccc}
\hline
$\tau$ & max $r_\textrm{dd}$ & $\textrm{size}(\mat{A}_\textrm{ff})$ & \% of unknowns \\
\hline
0.1 & 0.100 & 418  & 11.73 \\
0.2 & 0.198 & 524  & 14.71 \\
0.3 & 0.298 & 599  & 16.82 \\
0.4 & 0.398 & 709  & 19.90 \\
0.5 & 0.500 & 744  & 20.89 \\
0.6 & 0.598 & 815  & 22.88 \\
0.7 & 0.698 & 874  & 24.54 \\
0.8 & 0.797 & 948  & 26.61 \\
0.9 & 0.897 & 1012 & 28.41 \\
\hline
\end{tabular}
\end{table}
\subsection{Parallelism}
\label{sec:independent}
For AIRG, as noted in \secref{sec:cf_splitting}, $\tau$ controls how diagonally dominant, and hence how large, the resulting $\mat{A}_\textrm{ff}$ block is on each level, and we noted in \secref{sec:multilevel_ilu} that the size of $\mat{A}_\textrm{ff}$ sets the available parallelism. A large $\mat{A}_\textrm{ff}$ on each level, paired with an approximation of $\mat{A}_\textrm{ff}^{-1}$ that can be applied with only SpMVs, ensures there is sufficient fine-grained parallelism within the SpMVs on each level, and also gives fewer levels in the hierarchy. Many levels are undesirable, particularly on GPUs: each level requires a synchronisation and kernel launches on ever smaller coarse-grid systems, and the Schur complement coarse-grid matrices tend to become denser (there is fill-in even with dropping), giving more communication than work when distributed on lower levels. There is, however, no guarantee for a general linear system that a large diagonally dominant $\mat{A}_\textrm{ff}$ can be extracted on every level; for a dense $\mat{L}$ or $\mat{U}$ with entries all of similar magnitude the hierarchy can degenerate into a single-unknown Gaussian elimination with no parallelism, so reduction multigrids are best suited to sparse factors.

Increasing $\tau$ is one mechanism for reducing the number of levels, with another being the automatic truncation of the multigrid hierarchy we introduced in \cite{Dargaville2026}, which removes the coarsest levels, namely those with the least parallelism. On each level during the setup we generate a high-order GMRES polynomial as a tentative coarse-grid solver and test whether it can solve the coarse-grid problem to a loose relative tolerance against a random right-hand side; if it can, we truncate the hierarchy on that level and use the polynomial as the coarse-grid solver. We are able to truncate in this way precisely because AIRG only requires an approximate solve on each level. Truncation trades additional work for more parallelism (and a less dense communication pattern), typically resulting in a faster solve time.

Much as the number of level sets limits the parallelism of an exact triangular solve, the parallelism available in a multilevel hierarchy is problem-dependent. A level-set substitution must process its level sets in sequence, with the rows within each level set processed in parallel. Similarly, an AIRG hierarchy is an ordered sequence of levels, as the F-point smoothing on each level must be completed before proceeding to the next, while with an SpMV-based application of an $\mat{A}_\textrm{ff}^{-1}$ approximation all the F-points on a level are smoothed in parallel. We therefore compare the number of level sets of a factor with the number of levels in its AIRG hierarchy. These counts are not directly comparable: each level of a V-cycle applies a small, fixed number of SpMVs (with the assembled $q_{m-1}(\mat{A}_\textrm{ff})$, $\mat{A}_\textrm{fc}$, $\mat{R}$ and $\mat{P}$), each of which is fully parallel, so the number of sequential kernels in an AIRG application is a small multiple of the number of levels, compared with one per level set for a level-set solve. The actual cost of both is measured by the GPU timings in \secref{sec:axb_solves}. The ILU(0) of \textsc{bcsstk24} has 856 level sets in a natural ordering and 486 in \textsc{rcm}, an average of 4.2 and 7.3 rows per level set, respectively. For the $\mat{L}$ factor, AIRG with the DIAG\_DOM CF splitting and auto-truncation produces a hierarchy of only 24 levels, with 744 F-points (21\% of the unknowns) on the top grid.

These counts also behave differently under refinement. The number of levels in an untruncated multigrid hierarchy grows like $\mathcal{O}(\log n)$ with the number of unknowns $n$, while the number of level sets of a triangular factor can grow much faster; we give an example in \secref{sec:advection-diffusion}. A multilevel approach based on extracting diagonally dominant submatrices can therefore require far fewer sequential steps than a level-set approach in problems where the latter is close to serial (e.g., 1D upwind advection). We examine this below by reporting the number of level sets and of AIRG levels on several test problems.
\section{Results}
\label{sec:Results}
Unless otherwise noted, each factor $\mat{L}$ or $\mat{U}$ is solved with an undamped Richardson iteration preconditioned by the method we are testing. We monitor the unpreconditioned residual norm, use a random right-hand side, a zero initial guess, a relative tolerance of $10^{-6}$ and a maximum of 2000 iterations. One iteration is therefore one V-cycle of AIRG, one application of the GMRES polynomial or of the ISAI, or one Jacobi sweep. The cost of an iteration differs between the methods: a Jacobi sweep requires one SpMV with the factor, an ISAI iteration one SpMV with the factor and one with the approximate inverse (which has the same sparsity), a GMRES polynomial iteration seven SpMVs with the factor (we use 6\textsuperscript{th}-order polynomials), and an AIRG V-cycle several SpMVs on each level of its hierarchy. The parameters used for AIRG and the other iterative methods are fixed across all problems tested in this paper unless otherwise noted; the full set of solver and preconditioner options for PETSc and PFLARE is given in \cref{sec:parameters} of the supplementary materials. Given that iterative methods in non-normal systems can appear to diverge before converging, we must take care when interpreting convergence results. We therefore define divergence as occurring when a relative residual norm grows by more than $10^4$.

We consider two standard multigrid metrics when discussing AIRG: the grid complexity, defined as the total number of unknowns across all levels divided by the number on the finest level, and the storage complexity, defined as the total number of nonzeros in the operators that must be stored to apply a multigrid cycle divided by the number of nonzeros in the finest-level matrix (a proxy for the memory required by the hierarchy). As we only use F-point smoothing, the full operator on each level is not required during the solve; the stored operators on each level are $\mat{A}_\textrm{ff}$, $\mat{A}_\textrm{fc}$, $\hat{\mat{A}}_\textrm{ff}^{-1}$, $\mat{R}$ and $\mat{P}$, along with the coarsest matrix, so the storage complexity can be below one.

We begin by testing on three families of matrices directly from discretisations of PDEs often used as examples in the ILU literature, focusing on the achieved convergence. We sweep the grid size and parameter values; these are listed in \cref{sec:helmholtz,sec:advection-diffusion,sec:aniso_diffusion}. In all three cases the PDE is defined for a scalar unknown field $u(x,y)$ with a given source term $f(x,y)$ on the unit square $\Omega = (0,1)^2$, with Dirichlet boundary conditions. We assemble the discrete operator only on the $N \times N$ interior nodes of a uniform grid with spacing $h = 1/(N+1)$. The unknowns are assembled with the natural lexicographic ordering $\ell = j N + i$ for $i, j \in \{0, \ldots, N-1\}$, where $i$ indexes the $x$-direction and $j$ the $y$-direction.
\subsection{Indefinite Helmholtz}
\label{sec:helmholtz}
The shifted Helmholtz problem is
\begin{equation}
- \ppfrac{u}{x} - \ppfrac{u}{y} - k^2 u + \eta u = f \quad \text{on } \Omega,
\label{eq:helmholtz}
\end{equation}
where $k$ is the (real) wavenumber and $\eta \ge 0$ is a real shift. We discretise with a second-order 5-point central-difference stencil. The matrix is symmetric but indefinite, becoming more indefinite as $\eta \to 0$. Rather than fix $k$ in absolute terms, we parameterise the problem by the number of points per wavelength, $\text{ppw}$, and a shift percentage $s$, with
\begin{equation}
k = \frac{2 \pi}{\text{ppw} \cdot h}, \qquad \eta = \frac{s}{100} \, k^2.
\label{eq:helmholtz_params}
\end{equation}
Fixing $\text{ppw}$ across refinements means the wavelength always spans the same number of cells, so the problem is equally well resolved at every $N$; without this rescaling, refining the grid would make the problem easier (more points per wavelength). With $\text{ppw}$ fixed, grid refinement isolates the effect of the discrete problem size on the solver.

We sweep the grid size $N \in \{100, 200, 450, 800\}$ for $\text{ppw}=10$ with $s \in \{0, 5, 50\}$, and for $\text{ppw}=4$ with $s=0$ only, the fully indefinite case, which is the hardest problem in this family. For the ILU(0) factorisation every method converges robustly across the entire sweep, so we don't include its table. AIRG takes the fewest iterations, between 1 and 3, independent of $N$, the wavenumber and the shift, while the iteration counts of the other methods are also essentially flat under refinement (with Jacobi the largest, at up to 37 iterations).

The contrast appears for the ILU(1) factorisation, for which \tref{tab:results_helmholtz_ilu1} gives the iteration count required to solve the $\mat{L}$ and $\mat{U}$ factors. For $\text{ppw}=10$ all methods again converge with roughly constant iteration counts under refinement. For the fully indefinite, under-resolved $\text{ppw}=4$, $s=0$ case, however, the additional fill produces $\mat{L}$ and $\mat{U}$ factors that are much harder to solve iteratively, and the simpler methods degrade as the grid is refined. Jacobi diverges at every grid size. The GMRES polynomial converges only on the coarsest grid, and then only after more than 1700 iterations, hitting the maximum iteration count for every $N \ge 200$. ISAI is fairly resilient, converging for $N = 100$ and $N = 200$ before diverging for $N \ge 450$. AIRG, by contrast, converges in just 3 iterations at every grid size. This is the only Helmholtz case in which the simpler iterative methods fail, with AIRG the only method that remains robust throughout.
\begin{table}[ht]
\centering
\caption{For the indefinite Helmholtz problem \eref{eq:helmholtz} the number of iterations to solve the $\mat{L}$ and $\mat{U}$ factors generated by an ILU(1) factorisation to a relative tolerance of $10^{-6}$, as the grid is refined for each combination of points per wavelength ppw and shift $s$. ``MAX'' denotes hitting the maximum iteration count (2000) and ``DIV'' indicates the solve diverged.}
\label{tab:results_helmholtz_ilu1}
\small
\newcolumntype{G}{>{\centering\arraybackslash}p{2.5em}}
\begin{tabular}{c c r l G@{\hskip 0.6em}G!{\color{light}\vrule}G@{\hskip 0.6em}G!{\color{light}\vrule}G@{\hskip 0.6em}G!{\color{light}\vrule}G@{\hskip 0.6em}G}
\toprule
ppw & $s$ (\%) & $N$ & ordering & \multicolumn{2}{c}{AIRG} & \multicolumn{2}{!{\color{light}\vrule}c}{GMRES poly} & \multicolumn{2}{!{\color{light}\vrule}c}{ISAI} & \multicolumn{2}{!{\color{light}\vrule}c}{Jacobi} \\
\cmidrule(lr){5-6}\cmidrule(lr){7-8}\cmidrule(lr){9-10}\cmidrule(lr){11-12}
 &  &  &  & L & U & L & U & L & U & L & U \\
\midrule
\multirow{4}{*}{4} & \multirow{4}{*}{0} & 100 & \texttt{qmd} & 3 & 3 & 1754 & 1901 & 18 & 21 & \cellcolor{light}DIV & \cellcolor{light}DIV \\
 &  & 200 & \texttt{qmd} & 3 & 3 & \cellcolor{light}MAX & \cellcolor{light}MAX & 23 & 24 & \cellcolor{light}DIV & \cellcolor{light}DIV \\
 &  & 450 & \texttt{qmd} & 3 & 3 & \cellcolor{light}MAX & \cellcolor{light}MAX & \cellcolor{light}DIV & \cellcolor{light}DIV & \cellcolor{light}DIV & \cellcolor{light}DIV \\
 &  & 800 & \texttt{qmd} & 3 & 3 & \cellcolor{light}MAX & \cellcolor{light}MAX & \cellcolor{light}DIV & \cellcolor{light}DIV & \cellcolor{light}DIV & \cellcolor{light}DIV \\
\midrule
\multirow{4}{*}{10} & \multirow{4}{*}{0} & 100 & \texttt{colamd} & 3 & 3 & 5 & 5 & 16 & 15 & 31 & 30 \\
 &  & 200 & \texttt{colamd} & 3 & 3 & 7 & 6 & 19 & 17 & 37 & 35 \\
 &  & 450 & \texttt{mmd\_ata} & 3 & 3 & 8 & 8 & 24 & 22 & 48 & 45 \\
 &  & 800 & \texttt{nd} & 3 & 2 & 5 & 5 & 15 & 14 & 30 & 29 \\
\midrule
\multirow{4}{*}{10} & \multirow{4}{*}{5} & 100 & \texttt{nd} & 3 & 2 & 4 & 4 & 12 & 11 & 24 & 23 \\
 &  & 200 & \texttt{colamd} & 3 & 3 & 6 & 6 & 17 & 16 & 35 & 33 \\
 &  & 450 & \texttt{mmd\_ata} & 3 & 3 & 7 & 7 & 20 & 20 & 43 & 42 \\
 &  & 800 & \texttt{mmd\_ata} & 3 & 3 & 7 & 7 & 21 & 19 & 42 & 40 \\
\midrule
\multirow{4}{*}{10} & \multirow{4}{*}{50} & 100 & \texttt{rcm} & 3 & 3 & 13 & 13 & 35 & 35 & 73 & 73 \\
 &  & 200 & \texttt{rcm} & 3 & 3 & 13 & 13 & 35 & 35 & 74 & 74 \\
 &  & 450 & \texttt{natural} & 3 & 3 & 13 & 13 & 36 & 36 & 75 & 75 \\
 &  & 800 & \texttt{natural} & 3 & 3 & 13 & 13 & 36 & 36 & 75 & 75 \\
\bottomrule
\end{tabular}
\end{table}

\subsection{Advection-diffusion}
\label{sec:advection-diffusion}
We follow \cite{Chow2015} and consider the advection-diffusion problem
\begin{equation}
- \ppfrac{u}{x} - \ppfrac{u}{y} + \beta_\text{eff} \pfrac{(e^{xy} u)}{x} + \beta_\text{eff} \pfrac{(e^{-xy} u)}{y} = f \quad \text{on } \Omega,
\label{eq:advdiff}
\end{equation}
where $\beta_\text{eff} \ge 0$ is the advection coefficient that controls the strength of the convective term. The advection field $\mat{v}(x, y) = (e^{xy}, e^{-xy})$ is variable across the domain. We discretise both the diffusion and the advection terms with second-order central differences, giving a 5-point stencil. Central differencing of the advection term is used deliberately: for large $\beta$ it produces strongly nonsymmetric matrices far from diagonally dominant, whose ILU factors can be unstable, a standard hard test for ILU preconditioners \cite{Chow2015}.

To isolate the effect of grid refinement, we hold the cell P\'eclet number constant under refinement by setting
\begin{equation}
\beta_\text{eff} = \beta \cdot \frac{h_\text{ref}}{h},
\label{eq:beta_eff}
\end{equation}
where $h_\text{ref} = 1/(N_\text{ref}+1)$ is a fixed reference spacing with $N_\text{ref} = 450$ (the reference grid in \cite{Chow2015}). We sweep the grid size $N \in \{100, 200, 450, 800\}$ for each advection strength $\beta \in \{1500, 3000, 10000\}$; the values used in \cite{Chow2015} correspond to $\beta \in \{1500, 3000\}$ with $N=450$. \tref{tab:results_advdiff_ilu1} gives the number of iterations required to solve the $\mat{L}$ and $\mat{U}$ factors for the ILU(1) factorisation; we summarise the ILU(0) results in the text rather than including a table.

For the ILU(0) factorisation AIRG is robust and largely insensitive to refinement, converging in between 1 and 6 iterations across the whole sweep, with every method converging for every $\beta$ and $N$. The other methods, by contrast, are not scalable under refinement: for $\beta = 3000$ the Jacobi iteration count on the $\mat{L}$/$\mat{U}$ factors grows from $79/36$ at $N = 100$ to $144/206$ at $N = 800$, and ISAI similarly grows from $40/18$ to $72/103$, whereas AIRG stays between 3 and 6 iterations.

The ILU(1) factorisation (\tref{tab:results_advdiff_ilu1}) is far more challenging. The additional fill makes the $\mat{L}$ factor in particular very hard to solve: Jacobi diverges on the $\mat{L}$ factor for every $\beta$ and $N$ (and on both factors for $\beta = 10000$ and for $\beta = 3000$, $N = 800$), and ISAI diverges on $\mat{L}$ across almost the entire sweep. GMRES polynomials avoid divergence in all but one case (noted below) but are not scalable, with the $\mat{L}$ iteration count growing steeply under refinement, from $64$ to $327$ for $\beta = 1500$ and from $78$ to $512$ for $\beta = 10000$. AIRG is the only method that solves both factors with a low, near-constant iteration count throughout, taking 2--7 iterations for $\beta = 1500$ and $\beta = 3000$ and growing only mildly (from 7 to 17) for the hardest $\beta = 10000$ case; the lone exception is $\beta = 3000$, $N = 450$, where both AIRG and GMRES polynomials diverge on $\mat{L}$ in the chosen ordering. The \textsc{rcm} ILU(1) factor $\mat{L}$ at this size is far less stable than at $N = 200$ ($\|\mat{L}^{-1}\mat{b}\|/\|\mat{b}\|$ reaches $10^6$ for a random $\mat{b}$, compared with 180), and the residual of both methods initially grows by more than our divergence threshold of $10^4$. With the divergence check disabled, AIRG converges in 12 iterations and the GMRES polynomial in 272. We nonetheless count this as a failure: a single V-cycle (or a few iterations of the other methods) that amplifies the residual is of little use as a fixed-work inner solve for $\mat{A}\mat{x}=\mat{b}$ (\secref{sec:axb_solves}).

The iteration counts above are reported for the ordering that is best for the exact factorisation, chosen afresh at each $N$ (as noted in \secref{sec:ilu}); advection-dominated problems are well known to be sensitive to ordering due to long-range transport, and this choice can change the difficulty of the iterative factor solves considerably. For $\beta = 3000$ the ordering changes from \textsc{rcm} to \textsc{1wd} at $N = 800$, as the \textsc{rcm} ILU(1)-preconditioned $\mat{A}\mat{x}=\mat{b}$ solve diverges; this, rather than the grid size, is responsible for the drop in the GMRES polynomial iteration count on $\mat{L}$ from 127 at $N = 200$ to 10 at $N = 800$. The \textsc{1wd} factors are considerably easier for every iterative method, though worse as a preconditioner when solved exactly (the exact ILU(1)-preconditioned solve takes 36 to 350 iterations, compared with 18 to 119 for \textsc{rcm} up to $N = 450$). With \textsc{1wd} fixed at every $N$, AIRG solves the factors in 3 and 5 iterations independent of $N$, whereas the GMRES polynomial grows from 7 to 10 iterations on $\mat{L}$, ISAI from 48 to 102 on $\mat{U}$, and Jacobi diverges on both factors from $N = 450$. Conversely, for $\beta = 1500$, where the chosen ordering at $N = 450$ and $800$ is \textsc{natural}, fixing \textsc{rcm} instead gives ILU(1) factors $\mat{L}$ on which AIRG and the GMRES polynomial also diverge at those sizes.

\tref{tab:results_advdiff_ordering} isolates the effect of the ordering for the $\beta = 10000$, $N = 100$, ILU(0) problem. Of the nine orderings we tested, only three (\textsc{qmd}, \textsc{nd} and \textsc{colamd}) produce an ILU(0) factorisation for which the $\mat{A}\mat{x}=\mat{b}$ solve with exact triangular solves converges, so we restrict the comparison to these. Across all three, AIRG solves both factors in 1 to 3 iterations, whereas the iteration counts on $\mat{L}$/$\mat{U}$ of the GMRES polynomials range from $1/1$ under \textsc{nd} to $37/16$ under \textsc{colamd}, those of ISAI from $3/4$ to $58/45$, and those of Jacobi from $6/6$ to $115/90$. The choice of ordering, made for the exact factorisation, increases the difficulty of the factor solves by more than an order of magnitude for the single-level methods but leaves AIRG largely unaffected for these ILU(0) factors.

The number of AIRG levels also varies with the ordering, rising from 5 under \textsc{qmd} and \textsc{nd} to 7 under \textsc{colamd}; the easiest factors coincide with the shallowest hierarchies. It is likewise well known that the parallelism available to an exact triangular solve depends strongly on the ordering (e.g., see \cite{Liu1990, Chow2015, Chow2018} and the references therein), and the level-set counts here bear this out: a level-set substitution of these factors has 6 level sets under \textsc{nd} but 123 under \textsc{colamd}, compared with 7 AIRG levels. We also note that the parallelism available to a level-scheduled triangular solve degrades under refinement. With a natural ordering, the ILU(0) factors of this problem have $2N-1$ level sets, growing from 199 at $N = 100$ to 1599 at $N = 800$, whereas, as discussed in \secref{sec:independent}, the number of levels in a multigrid hierarchy grows only logarithmically.
\begin{table}[ht]
\centering
\caption{For the advection-diffusion problem \eref{eq:advdiff} the number of iterations to solve the $\mat{L}$ and $\mat{U}$ factors generated by an ILU(1) factorisation to a relative tolerance of $10^{-6}$, as the grid is refined for each advection strength $\beta$. ``MAX'' denotes hitting the maximum iteration count (2000) and ``DIV'' indicates the solve diverged.}
\label{tab:results_advdiff_ilu1}
\small
\newcolumntype{G}{>{\centering\arraybackslash}p{2.5em}}
\begin{tabular}{c r l G@{\hskip 0.6em}G!{\color{light}\vrule}G@{\hskip 0.6em}G!{\color{light}\vrule}G@{\hskip 0.6em}G!{\color{light}\vrule}G@{\hskip 0.6em}G}
\toprule
$\beta$ & $N$ & ordering & \multicolumn{2}{c}{AIRG} & \multicolumn{2}{!{\color{light}\vrule}c}{GMRES poly} & \multicolumn{2}{!{\color{light}\vrule}c}{ISAI} & \multicolumn{2}{!{\color{light}\vrule}c}{Jacobi} \\
\cmidrule(lr){4-5}\cmidrule(lr){6-7}\cmidrule(lr){8-9}\cmidrule(lr){10-11}
 &  &  & L & U & L & U & L & U & L & U \\
\midrule
\multirow{4}{*}{1500} & 100 & \texttt{rcm} & 5 & 3 & 64 & 3 & 103 & 10 & \cellcolor{light}DIV & 19 \\
 & 200 & \texttt{rcm} & 6 & 3 & 119 & 3 & \cellcolor{light}DIV & 10 & \cellcolor{light}DIV & 20 \\
 & 450 & \texttt{natural} & 5 & 2 & 189 & 3 & \cellcolor{light}DIV & 8 & \cellcolor{light}DIV & 18 \\
 & 800 & \texttt{natural} & 7 & 2 & 327 & 3 & \cellcolor{light}DIV & 8 & \cellcolor{light}DIV & 18 \\
\midrule
\multirow{4}{*}{3000} & 100 & \texttt{rcm} & 5 & 4 & 69 & 6 & \cellcolor{light}DIV & 27 & \cellcolor{light}DIV & 80 \\
 & 200 & \texttt{rcm} & 7 & 4 & 127 & 6 & \cellcolor{light}DIV & 28 & \cellcolor{light}DIV & 83 \\
 & 450 & \texttt{rcm} & \cellcolor{light}DIV & 4 & \cellcolor{light}DIV & 6 & \cellcolor{light}DIV & 28 & \cellcolor{light}DIV & 85 \\
 & 800 & \texttt{1wd} & 3 & 5 & 10 & 9 & 42 & 102 & \cellcolor{light}DIV & \cellcolor{light}DIV \\
\midrule
\multirow{4}{*}{10000} & 100 & \texttt{rcm} & 7 & 7 & 78 & 20 & \cellcolor{light}DIV & \cellcolor{light}DIV & \cellcolor{light}DIV & \cellcolor{light}DIV \\
 & 200 & \texttt{rcm} & 9 & 7 & 142 & 24 & \cellcolor{light}DIV & \cellcolor{light}DIV & \cellcolor{light}DIV & \cellcolor{light}DIV \\
 & 450 & \texttt{rcm} & 12 & 7 & 298 & 26 & \cellcolor{light}DIV & \cellcolor{light}DIV & \cellcolor{light}DIV & \cellcolor{light}DIV \\
 & 800 & \texttt{rcm} & 17 & 8 & 512 & 27 & \cellcolor{light}DIV & \cellcolor{light}DIV & \cellcolor{light}DIV & \cellcolor{light}DIV \\
\bottomrule
\end{tabular}
\end{table}

\begin{table}[ht]
\centering
\caption{The effect of the unknown ordering on the ILU(0) factors of the advection-diffusion problem \eref{eq:advdiff} with $\beta = 10000$, $N = 100$, for the three orderings (out of the nine tested) for which the exact ILU(0)-preconditioned $\mat{A}\mat{x}=\mat{b}$ solve converges.}
\label{tab:results_advdiff_ordering}
\small
\newcolumntype{G}{>{\centering\arraybackslash}p{2.1em}}
\begin{tabular}{l G@{\hskip 0.5em}G!{\color{light}\vrule}G@{\hskip 0.5em}G!{\color{light}\vrule}G@{\hskip 0.5em}G!{\color{light}\vrule}G@{\hskip 0.5em}G!{\color{light}\vrule}G@{\hskip 0.5em}G!{\color{light}\vrule}G}
\toprule
ordering & \multicolumn{8}{c}{iterations} & \multicolumn{3}{!{\color{light}\vrule}c}{levels} \\
\cmidrule(lr){2-9}\cmidrule(lr){10-12}
 & \multicolumn{2}{c}{AIRG} & \multicolumn{2}{!{\color{light}\vrule}c}{GMRES poly} & \multicolumn{2}{!{\color{light}\vrule}c}{ISAI} & \multicolumn{2}{!{\color{light}\vrule}c}{Jacobi} & \multicolumn{2}{!{\color{light}\vrule}c}{AIRG} & \multicolumn{1}{!{\color{light}\vrule}c}{level sets} \\
\cmidrule(lr){2-3}\cmidrule(lr){4-5}\cmidrule(lr){6-7}\cmidrule(lr){8-9}\cmidrule(lr){10-11}\cmidrule(lr){12-12}
 & L & U & L & U & L & U & L & U & L & U & \\
\midrule
\texttt{qmd} & 1 & 1 & 2 & 2 & 4 & 4 & 8 & 8 & 5 & 5 & 8 \\
\texttt{nd} & 2 & 2 & 1 & 1 & 3 & 4 & 6 & 6 & 5 & 5 & 6 \\
\texttt{colamd} & 3 & 3 & 37 & 16 & 58 & 45 & 115 & 90 & 7 & 6 & 123 \\
\bottomrule
\end{tabular}
\end{table}
\subsection{Anisotropic diffusion}
\label{sec:aniso_diffusion}
The anisotropic diffusion problem is
\begin{equation}
- \varepsilon \, \ppfrac{u}{x} - \ppfrac{u}{y} = f \quad \text{on } \Omega,
\label{eq:aniso}
\end{equation}
where $\varepsilon \in (0, 1]$ is the anisotropy ratio; $y$ is the strong direction and $x$ the weak direction. We use a second-order central-difference 5-point stencil, which gives a symmetric positive definite matrix, and we set $\varepsilon = 10^{-6}$ to produce very anisotropic diffusion. We sweep the grid size $N \in \{100, 200, 450, 800\}$ at this fixed anisotropy ratio. \tref{tab:results_aniso_ilu0} gives the number of iterations required to solve the $\mat{L}$ and $\mat{U}$ factors for the ILU(0) factorisation; the ILU(1) factorisation gives the same result, with every count within one of its ILU(0) value.
\begin{table}[ht]
\centering
\caption{For the anisotropic diffusion problem \eref{eq:aniso} with $\varepsilon = 10^{-6}$ the number of iterations to solve the $\mat{L}$ and $\mat{U}$ factors generated by an ILU(0) factorisation to a relative tolerance of $10^{-6}$, as the grid is refined (the winning ordering is \texttt{natural} at every $N$). ``MAX'' denotes hitting the maximum iteration count (2000) and ``DIV'' indicates the solve diverged.}
\label{tab:results_aniso_ilu0}
\small
\newcolumntype{G}{>{\centering\arraybackslash}p{2.5em}}
\begin{tabular}{r G@{\hskip 0.6em}G!{\color{light}\vrule}G@{\hskip 0.6em}G!{\color{light}\vrule}G@{\hskip 0.6em}G!{\color{light}\vrule}G@{\hskip 0.6em}G}
\toprule
$N$ & \multicolumn{2}{c}{AIRG} & \multicolumn{2}{!{\color{light}\vrule}c}{GMRES poly} & \multicolumn{2}{!{\color{light}\vrule}c}{ISAI} & \multicolumn{2}{!{\color{light}\vrule}c}{Jacobi} \\
\cmidrule(lr){2-3}\cmidrule(lr){4-5}\cmidrule(lr){6-7}\cmidrule(lr){8-9}
 & L & U & L & U & L & U & L & U \\
\midrule
100 & 2 & 2 & 40 & 39 & 50 & 50 & 100 & 100 \\
200 & 2 & 2 & 72 & 72 & 100 & 100 & 200 & 200 \\
450 & 2 & 2 & 150 & 149 & 225 & 225 & 450 & 450 \\
800 & 2 & 2 & 256 & 255 & 400 & 400 & 800 & 800 \\
\bottomrule
\end{tabular}
\end{table}

This problem shows the largest difference between the methods under refinement. For both the ILU(0) and ILU(1) factorisations, AIRG solves both factors in 2 iterations at every grid size, independent of $N$. Every other method, by contrast, has an iteration count that grows with the grid size: linearly for Jacobi and ISAI, and slightly more slowly for the GMRES polynomials. For ILU(0), the Jacobi count on each factor grows from 100 at $N=100$ to 800 at $N=800$, ISAI from 50 to 400 and the GMRES polynomial from 40 to 256. The strong-direction coupling produces long dependency chains in the factors, exactly the structure that makes a single-level iterative solve scale poorly under refinement, whereas the AIRG hierarchy reduces these to a constant number of iterations.

\subsection{SuiteSparse problems}
\label{sec:SuiteSparse}
We now test the iterative methods on 91 symmetric and nonsymmetric matrices drawn from the SuiteSparse Matrix Collection \cite{Davis2011} and the Matrix Market \cite{Boisvert1997}; they are listed, with their sizes and application areas, in \cref{tab:matrices_sym,tab:matrices_asym} of the supplementary materials. Many of these matrices were selected as they are tested in related papers, namely \cite{Chow2015, Liu2016a, Anzt2015, Anzt2018, Anzt2018a, Chow2018}. For every one of these matrices, GMRES(30) preconditioned by the ILU factorisation with exact triangular solves does not diverge, although it converges within 2000 iterations on only 66 of them for ILU(0) and 71 for ILU(1). The $\mat{L}$ and $\mat{U}$ solves in \secref{sec:lu_solves} are reported for all of the matrices, and the $\mat{A}\mat{x}=\mat{b}$ timings in \secref{sec:axb_solves} for those that converge. The results in this section were run on LUMI-G (CSC, Finland). We used one MPI rank, which on the AMD MI250X GPUs of LUMI-G corresponds to using one Graphics Compute Die (GCD, of which there are two per AMD MI250X) with 64 GB of high-bandwidth memory available per GCD. Using one MPI rank allows us to make easy comparisons against an exact triangular solve, namely the hipSPARSE triangular solve available through PETSc (\cref{sec:parameters} of the supplementary materials).

We configure PETSc with 64-bit indices in order to get timings representative of those at scale. For the iterative methods we use Kokkos Kernels SpMVs, rather than the hipSPARSE SpMVs; comparisons show the hipSPARSE SpMVs are typically faster on larger problems (up to 1.8$\times$), but slower on smaller ones.

\subsubsection{\texorpdfstring{$\mat{L}$ and $\mat{U}$}{L and U} factor solves}
\label{sec:lu_solves}
We begin by solving the $\mat{L}$ and $\mat{U}$ factors. The full per-matrix iteration counts for the symmetric and nonsymmetric matrices with ILU(0) and ILU(1) are given in \cref{tab:results_sym_ilu0,tab:results_sym_ilu1,tab:results_asym_ilu0,tab:results_asym_ilu1} of the supplementary materials.

For both the ILU(0) and ILU(1) factorisations, AIRG is very robust. We categorise failure as hitting the maximum number of iterations (2000), diverging on either of the $\mat{L}$ or $\mat{U}$ factors, or not finishing, either due to time (we allow a maximum time per matrix of 600 seconds) or memory limitations. For ILU(0), AIRG fails on only 2 of the 91 matrices. On \textsc{crankseg\_2} (\#74) it did not finish, as the run exceeded the time limit during the setup of the $\mat{U}$ factor (the setup of the $\mat{L}$ factor alone took 459 s), and it diverged on \textsc{bcsstk36} (\#63). We note that all other iterative methods also fail on that problem. ISAI fails on 3 matrices, all by diverging, with GMRES polynomials failing on 11 and Jacobi the least robust, failing on 21. For ILU(1), two matrices did not finish with AIRG, with no outright failures: the setup exceeded the time limit on both \textsc{crankseg} matrices. ISAI diverged on three, with the GMRES polynomial failing on 12 and Jacobi failing on 22 of the test matrices.

For problems where both the $\mat{L}$ and $\mat{U}$ factors could be solved, we take the larger of the $\mat{L}$ and $\mat{U}$ iteration counts for each matrix. These are relatively low for AIRG with ILU(0), ranging from 1 to 9. The range is larger for ISAI at 1--68, and much larger for the GMRES polynomial at 1--1619, with Jacobi at 2--157. For ILU(1) the iteration counts are generally comparable, with AIRG at 1--10, ISAI at 1--79, GMRES polynomials at 1--365 and Jacobi at 2--144. The maximum iteration counts would suggest the GMRES polynomials are a poor method, but they are driven by two outliers, \textsc{s2rmt3m1} (\#36) for ILU(0) and \textsc{bcsstk17} (\#49) for ILU(1). The median of these counts over the matrices on which each method succeeds is more reflective of the typical behaviour. For ILU(0) it is 3 for AIRG, 13 for ISAI, 5 for the GMRES polynomial and 26 for Jacobi, and for ILU(1) it is 3 for AIRG, 18 for ISAI, 8 for the GMRES polynomial and 37 for Jacobi.

The GMRES polynomials often succeed in relatively few iterations on problems where Jacobi diverges, such as \textsc{gyro\_m} (\#48), \textsc{wathen100} (\#50), \textsc{wathen120} (\#54), \textsc{Pres\_Poisson} (\#61), \textsc{crankseg\_1} (\#73) and \textsc{crankseg\_2} (\#74) for ILU(0). The GMRES polynomials fail on more problems than ISAI, and largely on different ones: for both ILU(0) and ILU(1) the GMRES polynomial fails on 10 matrices on which ISAI succeeds, while it succeeds on two of the three matrices on which ISAI fails for ILU(0), \textsc{aft01} (\#21) and \textsc{nd3k} (\#66), and on \textsc{nd3k} for ILU(1).

The robustness of AIRG comes at a cost, namely the time required to set up the multigrid hierarchy and the memory required to store it. Both Jacobi and the GMRES polynomials require no additional storage, while the ISAI requires one extra matrix per factor, which in this work we restrict to the sparsity of the original factors. AIRG typically requires more storage for the multigrid hierarchy, but this is problem dependent.

We compute the grid and storage complexity of the AIRG hierarchies built on the $\mat{L}$ and $\mat{U}$ factors for each problem, taken as the maximum over the two factors. For ILU(0) the grid complexity ranges from 1.05 on \textsc{Chem97ZtZ} (\#3) to 26.93 on \textsc{s3rmq4m1} (\#41), with 62 of the 90 hierarchies below 3 (57 of the 89 for ILU(1)), reflecting how easily large diagonally dominant submatrices can be extracted in different problems (\secref{sec:comparison}). The storage complexity ranges from 1.47 on \textsc{Trefethen\_20000} (\#53) to 24.65 on \textsc{cfd1} (\#64) and exceeds 10 on 15 of the hierarchies (16 for ILU(1)), in particular for dense structural problems such as the cylindrical shell matrices. The ILU(1) ranges are 1.05--37.90 and 0.39--23.96, respectively; the storage complexity below one (again on \textsc{Trefethen\_20000}) arises because only the operators needed for F-point smoothing are stored on each level (\secref{sec:Results}). The median complexities are more representative: for ILU(0) the median grid complexity is 2.52 and the median storage complexity 5.03, with ILU(1) similar at 2.68 and 6.35.

Before moving on to the solve time, we compare, for the $\mat{L}$ factor of each ILU(0) test matrix, the number of AIRG levels against the number of level sets (\secref{sec:independent}). AIRG builds fewer levels than there are level sets in all but one of the test matrices; the single exception is \textsc{Chem97ZtZ} (\#3), where both have two.

The key observation is that the number of AIRG levels is more constant across the range of problems; AIRG builds a hierarchy of between 2 and 108 levels (a median of 5). The median number of levels of 5 here matches the level at which we start testing for auto-truncation. For many of the matrices tested, a GMRES polynomial can form a sufficiently good coarse-grid solver to allow such heavy truncation. Without auto-truncation the median number of levels is 13. It is also not just small matrices that benefit from heavy truncation. For example, \textsc{af\_shell3} (\#75) and \textsc{CoupCons3D} (\#90) both require 19 multigrid levels without truncation, but the hierarchy can be truncated at 5 levels without affecting convergence. As noted above, truncating a multigrid hierarchy increases the amount of work required in the solve, as the coarsest matrix is larger, but given our coarse-grid solve is a GMRES polynomial this work is very parallel. We see this reflected in the solve times, which decrease with auto-truncation, with a median speedup of between 1.4$\times$ and 1.5$\times$ across the $\mat{L}$ and $\mat{U}$ factors at both fill levels.

For a level-set approach, the number of level sets has a median of 529 and reaches 20133 on \textsc{Transport} (\#91), whose AIRG hierarchy has only 5 levels. In some problems the level-set approach is entirely serial, e.g., each of the 2000 level sets of the $\mat{L}$ factor of \textsc{Trefethen\_2000} (\#11, $n=2000$) is a single row, against 5 AIRG levels. The ILU(1) factors show similar behaviour, made more pronounced by the additional fill, which lengthens the level-set dependency chains (a median of 828 level sets, compared with a median of 5 AIRG levels), and the same is true of the setup costs and $\mat{A}\mat{x}=\mat{b}$ solve times discussed below. AIRG therefore requires far fewer sequential steps than an exact triangular solve, even on problems where the exact solve has low parallelism. The single-level methods (Jacobi, ISAI and the GMRES polynomials) have no level dependencies at all, as each of their iterations consists of fully parallel SpMVs, so their number of sequential steps is set by their iteration count.

\begin{figure}[ht]
\centering
\includegraphics[width=0.98\textwidth]{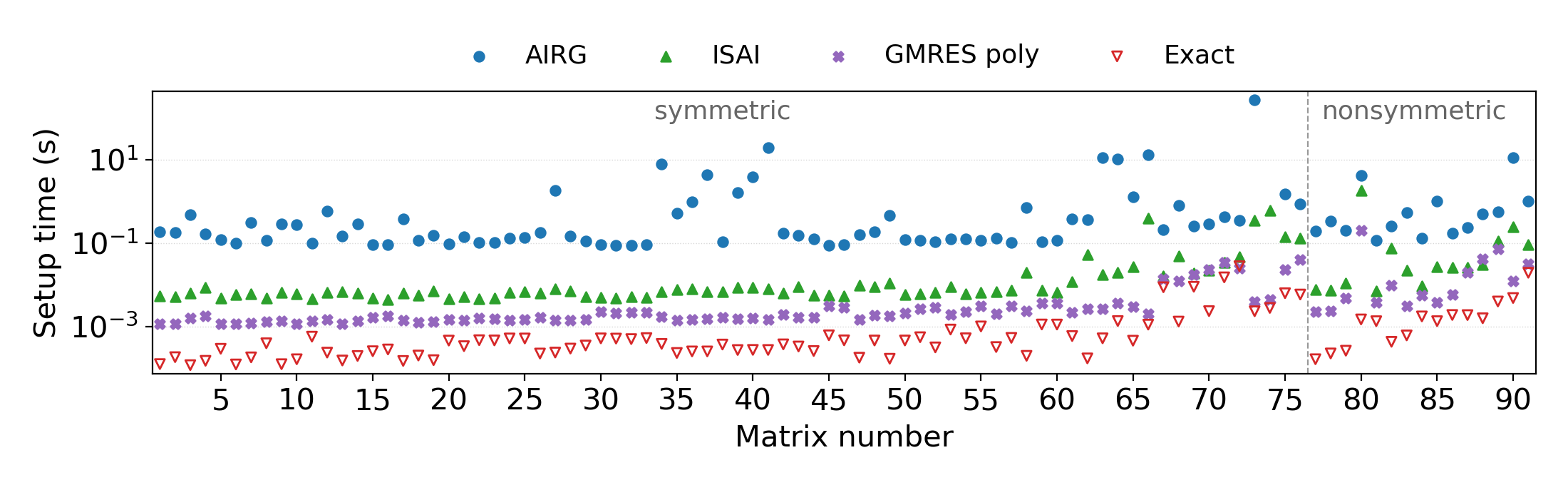}
\caption{Setup time comparison (smaller is better) between the GMRES polynomial, ISAI, AIRG and the analysis phase of the hipSPARSE exact triangular solve for the $\mat{L}$ and $\mat{U}$ factors generated by the ILU(0) factorisation of each test matrix, ordered as in the results tables (symmetric and nonsymmetric matrices separated by the dashed line). The results shown are taken as the maximum over the $\mat{L}$ and $\mat{U}$ factors. Jacobi requires no setup and is not shown.}
\label{fig:setup}
\end{figure}

All methods except Jacobi, however, cost time to set up, including the exact triangular solve, whose level-set analysis phase must run before the first solve. \fref{fig:setup} compares the setup cost of the different methods for ILU(0), again taken as the maximum over the $\mat{L}$ and $\mat{U}$ factors (the ILU(1) counterpart is given in \cref{sec:sm_figures} of the supplementary materials). The GMRES polynomial is the cheapest of the approximate methods to set up as it only requires generating the polynomial coefficients, whereas AIRG is the most expensive as it must build the multigrid hierarchy, with ISAI in between. The median setup cost for ILU(0) is 0.46 ms for the exact triangular solve, 2.0 ms for the GMRES polynomial, 7.2 ms for ISAI and 185 ms for AIRG, with the ILU(1) costs at 0.57 ms, 2.1 ms, 9.6 ms and 258 ms, respectively. The setup costs of AIRG also have a long tail: it takes more than 10 s on 6 of the problems for ILU(0), up to 268 s on \textsc{crankseg\_1} (\#73) (the $\mat{L}$ factor of \textsc{crankseg\_2} (\#74) alone takes 459 s), and on 9 for ILU(1), up to 56.9 s on \textsc{CoupCons3D} (\#90), excluding the two \textsc{crankseg} matrices, which did not finish. ISAI exceeds 10 s only for ILU(1), on \textsc{nd3k} (\#66, 44.5 s) and \textsc{crankseg\_2} (\#74), while the GMRES polynomial and the exact triangular solve never take more than 0.25 s. These setup times exclude the ILU factorisation itself and the ordering sweep of \secref{sec:ilu}, which are the same for every method.

\subsubsection{Solving \texorpdfstring{$\mat{A}\mat{x}=\mat{b}$}{Ax=b}}
\label{sec:axb_solves}
We now move to solving the $\mat{A}\mat{x}=\mat{b}$ system, using an outer GMRES(30) preconditioned by the ILU factorisation, with the inner application of $\mat{L}^{-1}$ and $\mat{U}^{-1}$ performed by iterative methods as detailed in the previous section. We use a random right-hand side and an initial guess of zero, and solve to a relative tolerance of $10^{-6}$, with a maximum number of outer iterations set to 2000. To ensure fair comparisons, we do not plot timing results for problems that diverged or reached the limit of 2000 outer iterations. The time of a solve that reaches the limit is just 2000 times the cost of one outer iteration; comparing such times would compare only the cost per iteration of each inner solve and ignore the number of iterations needed to converge, which \secref{sec:lu_solves} showed varies widely between methods and problems. As a baseline we show the timing results from using an exact inner triangular solve.

We know from \cite{Chow2015, Anzt2015, Chow2018} that a few iterations of an iterative method are sufficient in some problems to converge the inner solve well enough to not (overly) affect the outer iteration count. Choosing how many inner iterations to apply, however, is problem dependent; we could choose to solve each inner problem to a set tolerance but that would require computing a norm at each iteration (we would like to avoid dot products in the inner solve) and would require a flexible Krylov method for the outer solve (increasing the memory required). Instead, we choose a fixed number of iterations for the inner solve.

As AIRG solves either factor to a relative tolerance of $10^{-6}$ in a median of 3 iterations (and at most 10, \secref{sec:lu_solves}), we use only a single V-cycle for both $\mat{L}$ and $\mat{U}$ on the inner solve. This gives the same median number of outer iterations as an exact triangular solve for both ILU(0) and ILU(1). Picking the best number of inner iterations for the Jacobi, ISAI and GMRES polynomial methods is more difficult, as their iteration counts when solving $\mat{L}$ or $\mat{U}$ to $10^{-6}$ vary from a handful to hundreds, as noted in \secref{sec:lu_solves}. To determine a fair number, we ran each test matrix for both fill levels with each of these methods using a fixed number of inner iterations from $\{1, 2, 3, 5, 8, 13, 21, 34, 55, 89, 144\}$. For each problem and method we use the number of inner iterations that resulted in the smallest runtime while still converging. All the solve times reported in this section are from warmed-up runs, so one-time GPU initialisation costs (which can be comparable to the solve time of the smaller problems) are excluded.

\tref{tab:axb_summary} summarises the results, with the solve time for each test matrix given in \cref{fig:axb_timing_ilu0,fig:axb_timing_ilu1} of the supplementary materials. AIRG is faster than the hipSPARSE exact triangular solve on 42 of the 65 problems where both converge for ILU(0) and 52 of the 69 for ILU(1), with median speedups of 1.3$\times$ and 2.2$\times$, respectively. The larger setup cost of AIRG is not repaid immediately, however: on the problems where AIRG is faster, the extra setup time for both factors (relative to the analysis phase of the exact triangular solve) is recovered after a median of 11.3 solves of $\mat{A}\mat{x}=\mat{b}$ for ILU(0) and 11.7 for ILU(1), and within 10 solves on 21 of the 42 and 25 of the 52 such problems, respectively. The tuned Jacobi, ISAI and GMRES polynomial inner solves are faster than the exact triangular solve on all but 12, 6 and 9, respectively, of the problems where both converge for ILU(0), and 9, 7 and 9 for ILU(1).

\begin{table}[ht]
\centering
\caption{Summary of the $\mat{A}\mat{x}=\mat{b}$ solves on the 91 SuiteSparse/Matrix Market problems. ``Conv.'' is the number of problems solved, ``Fastest'' the number on which the method had the smallest solve time, and ``Speedup'' the median speedup over the exact triangular solve on the problems where both converged, over all of them and over those with $n\geq10^5$. Jacobi, ISAI and the GMRES polynomial use the best number of inner iterations for each problem.}
\label{tab:axb_summary}
\small
\begin{tabular}{l rrrr rrrr}
\toprule
 & \multicolumn{4}{c}{ILU(0)} & \multicolumn{4}{c}{ILU(1)} \\
\cmidrule(lr){2-5}\cmidrule(lr){6-9}
 & & & \multicolumn{2}{c}{Speedup} & & & \multicolumn{2}{c}{Speedup} \\
\cmidrule(lr){4-5}\cmidrule(lr){8-9}
Method & Conv. & Fastest & All & $n\geq10^5$ & Conv. & Fastest & All & $n\geq10^5$ \\
\midrule
Exact & 66 & 4 & 1 & 1 & 71 & 4 & 1 & 1 \\
AIRG & 65 & 1 & 1.3 & 1.2 & 69 & 5 & 2.2 & 2.2 \\
Jacobi & 63 & 3 & 2.6 & 2.2 & 65 & 1 & 3.2 & 2.2 \\
GMRES poly & 67 & 9 & 3.1 & 2.4 & 71 & 11 & 5.0 & 2.0 \\
ISAI & 67 & 50 & 4.5 & 2.7 & 70 & 51 & 6.8 & 2.6 \\
\bottomrule
\end{tabular}
\end{table}

Their speedups are also larger than those of AIRG, with ISAI the fastest: it is the fastest of all the methods on 50 of the problems for ILU(0) and 51 for ILU(1), with median speedups over the exact solve of 4.5$\times$ and 6.8$\times$, respectively. Tuned in this way, ISAI and the GMRES polynomials also converge on slightly more problems than AIRG, while Jacobi converges on fewer. AIRG converges on every problem on which the exact triangular solve converges, except \textsc{crankseg\_2} (\#74) for ILU(0) and both \textsc{crankseg} matrices for ILU(1). ISAI and the GMRES polynomials converge on one more problem than the exact triangular solve with ILU(0), and with ILU(1) the GMRES polynomial converges on one problem on which the exact triangular solve does not: on \textsc{af\_shell3} (\#75) with ILU(0) and \textsc{apache2} (\#69) with ILU(1) the exact triangular solve stops just short of the tolerance at the limit of 2000 outer iterations, and the tuned methods converge there just inside the limit. The GMRES polynomial does not converge on \textsc{cbuckle} (\#58) with ILU(1), however, so overall it converges on the same number of problems as the exact triangular solve.

The difference between the methods, however, depends strongly on the problem size. The median speedup of AIRG over the exact triangular solve is lowest on the smallest problems: 1.1$\times$ for the problems with $n<20\,000$, against 2.0$\times$ for $20\,000\leq n<100\,000$ and 1.2$\times$ for $n\geq100\,000$ with ILU(0), and 1.6$\times$, 2.7$\times$ and 2.2$\times$, respectively, with ILU(1). The speedups of ISAI and the GMRES polynomials fall on the largest problems, most of all for ISAI, from 4.8$\times$ to 2.7$\times$ with ILU(0) and from 6.7$\times$ to 2.6$\times$ with ILU(1), while that of Jacobi changes little. Comparing ISAI directly with AIRG on the problems where both converge, ISAI is a median of 4.5$\times$ faster than AIRG for $n<20\,000$ but only 2.0$\times$ faster for $n\geq100\,000$ with ILU(0), and 4.0$\times$ and 1.2$\times$ faster with ILU(1); on the largest problems with ILU(1) all four iterative methods give similar median speedups, between 2.0$\times$ and 2.6$\times$. The same trend holds when the problems are split by their nnz. This is consistent with the solves of the smaller problems being limited by the latency of launching GPU kernels rather than by the work they perform: an AIRG V-cycle launches several kernels on every level of the hierarchy, whereas the other methods apply only a few sparse matrix-vector products, and on the larger problems this fixed cost is amortised over more work.

These timings, however, rely on tuning the number of inner iterations of the simpler methods for each problem, which required solving each problem with 11 different inner iteration counts. The best number is not known in advance and varies widely between problems, as the factor solves in \secref{sec:lu_solves} showed, and a poor choice can stop the outer solve from converging altogether.
\section{Conclusions}
In this paper we have compared four parallel iterative methods for solving the $\mat{L}$ and $\mat{U}$ triangular factors from ILU factorisations, namely a Jacobi method, ISAI, GMRES polynomials and the reduction multigrid AIRG, both on the factors alone and as the inner solve of an ILU preconditioner on a single GPU. All of the iterative methods compared in this work are available in the open-source PFLARE library, which can be installed directly through the PETSc configure with \texttt{--download-pflare}.

To improve the convergence of AIRG in difficult problems, we introduced the DIAG\_DOM CF splitting, which enforces a user-specified maximum diagonal dominance ratio on $\mat{A}_\textrm{ff}$ with a single parameter, $\tau$, and is built entirely from parallel independent set routines, making it well suited to GPUs. AIRG requires far fewer sequential steps than an exact triangular solve: its hierarchies have a median of 5 levels on the tested problems, each requiring only a small number of fully parallel SpMVs, compared with a median of 529 level sets for an exact triangular solve of the ILU(0) factors (828 for ILU(1)), and on some matrices an exact substitution is essentially serial. The automatic truncation of the AIRG hierarchy we introduced previously further reduces the number of levels.

As the inner solve of an ILU preconditioner, the iterative methods gave only modest speedups over the hipSPARSE exact triangular solve. One V-cycle of AIRG gave median speedups of 1.3$\times$ for ILU(0) and 2.2$\times$ for ILU(1), and where it was faster its larger setup cost took a median of 11.3 and 11.7 solves of $\mat{A}\mat{x}=\mat{b}$ to amortise, respectively. With the number of inner iterations tuned for each problem, the other methods were faster, with ISAI the fastest method on most problems and median speedups over the exact triangular solve of 4.5$\times$ for ILU(0) and 6.8$\times$ for ILU(1). This advantage comes mainly from the smaller problems, where kernel launch times can dominate: on the problems with $n\geq100\,000$ the speedups of the four methods are closer. For example, for ILU(1), all four iterative methods produced median speedups of between 2.0$\times$ and 2.6$\times$ compared with the exact triangular solve.

We note that these comparisons are on a single GPU; exact triangular solves often perform poorly with distributed-memory parallelism, which is why distributed ILU preconditioners commonly fall back to block Jacobi with local factors. In contrast, all four iterative methods tested in this work extend directly to distributed-memory parallelism (with GPU-aware MPI used in the setup and solve phases in PFLARE). In particular, AIRG has been shown to weak scale well up to 11k GPUs when solving triangular systems from discretisations of advection problems \cite{Dargaville2026}.

Overall, this work shows that each method has advantages and can provide speedups when compared to exact triangular solves. The Jacobi method is the simplest, with no setup and no additional storage, but is the least robust, failing on the most problems both when solving the factors and as an inner solve. The GMRES polynomials share these advantages but are more robust: they fail on roughly half as many problems when solving the factors, and as an inner solve they converge on every problem on which the Jacobi method converges (and on several more) with a similar or faster solve time, making them a potential replacement for a Jacobi method. ISAI gives the largest speedups as an inner solve and is nearly as robust as AIRG on the SuiteSparse/Matrix Market factors, at the cost of one extra stored matrix per factor and a larger setup than the Jacobi method or GMRES polynomials, although it diverges on most of the ILU(1) advection-diffusion problems. AIRG is the most robust method when solving the factors, with nearly flat iteration counts under refinement where those of every other method grow or diverge, and needs no per-problem tuning. Its drawbacks are memory use, at a median of five to six copies of the factor, and the largest setup times. Which method is best therefore depends on whether robustness, solve time or memory use matters most, and on whether the inner solve is tuned for the problems of interest.

\section*{Acknowledgments}
We acknowledge the EuroHPC Joint Undertaking for awarding this project access to the EuroHPC supercomputer LUMI, hosted by CSC (Finland) and the LUMI consortium through a EuroHPC Benchmark Access call.
\bibliographystyle{siamplain}
\bibliography{bib_library}

\clearpage
\makeatletter
\def\siampretitle{Supplementary Materials: }
\def\siamprelabel{SM}
\makeatother
\setcounter{section}{0}
\setcounter{figure}{0}
\setcounter{table}{0}
\setcounter{equation}{0}
\setcounter{page}{1}
\thispagestyle{plain}
\headers{Iterative methods for ILU solves}{S. Dargaville and C. C. Pain}
\begin{center}
{\bfseries SUPPLEMENTARY MATERIALS: A COMPARISON OF PARALLEL ITERATIVE METHODS FOR ILU TRIANGULAR SOLVES\par}
\vskip .075in
{\footnotesize STEVEN DARGAVILLE {\scriptsize AND} CHRISTOPHER C. PAIN\par}
\end{center}
\vskip .11in

\section{Solver parameters}
\label{sec:parameters}
All results use PFLARE commit \texttt{fef661d} and PETSc commit \texttt{6f69904}; these include optimisations to the GPU setup and solve of AIRG, the GMRES polynomials and ISAI made since \cite{Dargaville2026}.
The factor solves described at the start of \cref{sec:Results} of the article use the following command-line options with PETSc:
\begin{lstlisting}
-ksp_type richardson -ksp_norm_type unpreconditioned -ksp_rtol 1e-6 -ksp_atol 1e-50 -ksp_max_it 2000
\end{lstlisting}

When used as the inner solve of the ILU preconditioner for $\mat{A}\mat{x}=\mat{b}$ (\cref{sec:axb_solves} of the article), each factor is instead solved with a fixed number of iterations $k$ and without computing any norms or dot products: for $k=1$ (and always for the single V-cycle of AIRG) we use
\begin{lstlisting}
-ksp_type preonly
\end{lstlisting}
and for $k>1$
\begin{lstlisting}
-ksp_type richardson -ksp_norm_type none -ksp_max_it k
\end{lstlisting}
so that each inner solve is a fixed linear operator. The outer solve uses
\begin{lstlisting}
-ksp_type gmres -ksp_gmres_restart 30 -ksp_norm_type unpreconditioned -ksp_rtol 1e-6 -ksp_atol 1e-50 -ksp_max_it 2000
\end{lstlisting}
and the exact triangular solve baseline uses PETSc's ILU factorisation (inside block Jacobi with one block, as we use one MPI rank), with the matrix converted to PETSc's hipSPARSE format so that the triangular solves use hipSPARSE:
\begin{lstlisting}
-mat_type aijhipsparse -vec_type hip -pc_type bjacobi -sub_pc_type ilu -sub_pc_factor_mat_solver_type hipsparse
\end{lstlisting}
with \texttt{-sub\_pc\_factor\_levels 1} for ILU(1).

When using PFLARE we specify the matrix/vector types as the AIJKOKKOS format, since this ensures all the components of the setup and solve are run on the GPU. For AIRG we then also use:
\begin{lstlisting}
-pc_type air -pc_air_cf_splitting_type diag_dom -pc_air_strong_threshold 0.5 -pc_air_a_drop 1e-6 -pc_air_r_drop 1e-4 -pc_air_diag_scale_polys -pc_air_coarsest_diag_scale_polys -pc_air_coarsest_inverse_type newton -pc_air_coarsest_poly_order 10 -pc_air_coarsest_matrix_free_polys -pc_air_auto_truncate_start_level 5 -pc_air_auto_truncate_tol 1e-3
\end{lstlisting}
which uses the default 6\textsuperscript{th}-order GMRES polynomial ($m=7$) on all levels except the coarsest. With the DIAG\_DOM CF splitting, \texttt{-pc\_air\_strong\_threshold} sets both the strength threshold $\theta$ and the maximum diagonal dominance ratio $\tau$ of \cref{alg:pmisr_dd} of the article ($\theta=\tau$); the value $\tau=0.5$ used throughout is also the default. The last six options concern the coarse grid. Rather than coarsening until the problem is small enough to solve directly, from the fifth level onwards we test whether the coarse-grid solver alone can solve the current coarse problem to a relative tolerance of $10^{-3}$, and if so we truncate the hierarchy there. The coarse-grid solver is a 10\textsuperscript{th}-order GMRES polynomial in a Newton basis applied matrix-free, so it requires only matrix-vector products. This trades a higher cycle complexity for fewer levels, which is beneficial on GPUs where the lower levels of the hierarchy have little work but still incur kernel-launch costs.

The GMRES polynomial approximate inverse uses:
\begin{lstlisting}
-pc_type pflareinv -pc_pflareinv_type newton -pc_pflareinv_matrix_free
\end{lstlisting}
again with the default 6\textsuperscript{th}-order polynomial, with coefficients computed in a Newton basis and applied matrix-free.

The ISAI approximate inverse uses:
\begin{lstlisting}
-pc_type pflareinv -pc_pflareinv_type isai
\end{lstlisting}
which defaults to using the sparsity pattern of the input matrix.

The Jacobi iteration uses:
\begin{lstlisting}
-pc_type jacobi
\end{lstlisting}

When solving $\mat{A}\mat{x}=\mat{b}$ with an ILU preconditioner whose triangular factors are applied approximately (\secref{sec:axb_solves}), the factor solves above are instead applied as a fixed amount of work rather than to a tolerance, e.g., a single AIRG V-cycle or a fixed number of Jacobi iterations on the $\mat{L}$ and $\mat{U}$ factors.
\section{Test matrices}
\label{sec:sm_matrices}
\Cref{tab:matrices_sym,tab:matrices_asym} list the 91 SuiteSparse/Matrix Market test matrices used in \cref{sec:SuiteSparse} of the article, in the order used in the article's figures and in the results tables below (the symmetric matrices are numbered 1 to 76 and the nonsymmetric matrices 77 to 91, each group sorted by the number of nonzeros), with the size $n$ and number of nonzeros (nnz) of the full matrix and the application area recorded by the collections.
\par\begingroup\footnotesize
\begin{longtable}{rlrrl}
\caption{The symmetric SuiteSparse/Matrix Market test matrices (numbers 1 to 76), with the size $n$ and number of nonzeros (nnz) of the full matrix and the application area recorded by the collections. The numbering is that used throughout the article.}
\label{tab:matrices_sym} \\
\toprule
\# & Matrix & $n$ & nnz & Application area \\
\midrule
\endfirsthead
\toprule
\# & Matrix & $n$ & nnz & Application area \\
\midrule
\endhead
\midrule
\endfoot
\bottomrule
\endlastfoot
1 & \texttt{1138\_bus} & 1\,138 & 4\,054 & Power-system network \\
2 & \texttt{bcsstk19} & 817 & 6\,853 & Structural (part of a suspension bridge) \\
3 & \texttt{Chem97ZtZ} & 2\,541 & 7\,361 & Statistics (normal equations) \\
4 & \texttt{bcsstk08} & 1\,074 & 12\,960 & Structural (TV studio) \\
5 & \texttt{mhd3200b} & 3\,200 & 18\,316 & Magnetohydrodynamics \\
6 & \texttt{bcsstm12} & 1\,473 & 19\,659 & Structural (mass matrix) \\
7 & \texttt{bcsstk10} & 1\,086 & 22\,070 & Structural mechanics (buckling) \\
8 & \texttt{mhd4800b} & 4\,800 & 27\,520 & Magnetohydrodynamics \\
9 & \texttt{msc01050} & 1\,050 & 29\,156 & Structural mechanics \\
10 & \texttt{plbuckle} & 1\,282 & 30\,644 & Structural (buckling plate) \\
11 & \texttt{Trefethen\_2000} & 2\,000 & 41\,906 & Combinatorics (Trefethen) \\
12 & \texttt{bcsstk27} & 1\,224 & 56\,126 & Structural (buckling) \\
13 & \texttt{bcsstk14} & 1\,806 & 63\,454 & Structural (roof of Omni Coliseum) \\
14 & \texttt{nasa2146} & 2\,146 & 72\,250 & Structural mechanics \\
15 & \texttt{fv1} & 9\,604 & 85\,264 & 2D finite-volume mesh \\
16 & \texttt{fv2} & 9\,801 & 87\,025 & 2D finite-volume mesh \\
17 & \texttt{msc04515} & 4\,515 & 97\,707 & Structural mechanics \\
18 & \texttt{crystm01} & 4\,875 & 105\,339 & Materials (crystal mass matrix) \\
19 & \texttt{bcsstk15} & 3\,948 & 117\,816 & Structural (module of an offshore platform) \\
20 & \texttt{bodyy4} & 17\,546 & 121\,938 & Structural (FE mesh) \\
21 & \texttt{aft01} & 8\,205 & 125\,567 & Acoustics (vibroacoustics) \\
22 & \texttt{bodyy5} & 18\,589 & 129\,281 & Structural (FE mesh) \\
23 & \texttt{bodyy6} & 19\,366 & 134\,748 & Structural (FE mesh) \\
24 & \texttt{ted\_B} & 10\,605 & 144\,579 & Thermal model reduction \\
25 & \texttt{ted\_B\_unscaled} & 10\,605 & 144\,579 & Thermal model reduction (unscaled) \\
26 & \texttt{bcsstk18} & 11\,948 & 149\,090 & Structural (R.E. Ginna nuclear power station) \\
27 & \texttt{bcsstk24} & 3\,562 & 159\,910 & Structural (Calgary Olympic Saddledome arena) \\
28 & \texttt{Muu} & 7\,102 & 170\,134 & Structural (FEM mass matrix) \\
29 & \texttt{t2dah\_e} & 11\,445 & 176\,117 & Model reduction (micropyros thruster) \\
30 & \texttt{obstclae} & 40\,000 & 197\,608 & Obstacle problem (optimisation) \\
31 & \texttt{torsion1} & 40\,000 & 197\,608 & Torsion (optimisation) \\
32 & \texttt{jnlbrng1} & 40\,000 & 199\,200 & Journal-bearing (optimisation) \\
33 & \texttt{minsurfo} & 40\,806 & 203\,622 & Minimal-surface (optimisation) \\
34 & \texttt{s3rmt3m3} & 5\,357 & 207\,695 & Structural (cylindrical shell) \\
35 & \texttt{s1rmt3m1} & 5\,489 & 219\,521 & Structural (cylindrical shell) \\
36 & \texttt{s2rmt3m1} & 5\,489 & 219\,521 & Structural (cylindrical shell) \\
37 & \texttt{s3rmt3m1} & 5\,489 & 219\,521 & Structural (cylindrical shell) \\
38 & \texttt{Dubcova1} & 16\,129 & 253\,009 & 2D/3D PDE \\
39 & \texttt{s1rmq4m1} & 5\,489 & 281\,111 & Structural (cylindrical shell) \\
40 & \texttt{s2rmq4m1} & 5\,489 & 281\,111 & Structural (cylindrical shell) \\
41 & \texttt{s3rmq4m1} & 5\,489 & 281\,111 & Structural (cylindrical shell) \\
42 & \texttt{chipcool0} & 20\,082 & 281\,150 & Thermal model reduction \\
43 & \texttt{bcsstk16} & 4\,884 & 290\,378 & Structural (Corps of Engineers dam) \\
44 & \texttt{crystm02} & 13\,965 & 322\,905 & Materials (crystal mass matrix) \\
45 & \texttt{shallow\_water1} & 81\,920 & 327\,680 & Climate (shallow-water) \\
46 & \texttt{shallow\_water2} & 81\,920 & 327\,680 & Climate (shallow-water) \\
47 & \texttt{Kuu} & 7\,102 & 340\,200 & Structural (FEM stiffness) \\
48 & \texttt{gyro\_m} & 17\,361 & 340\,431 & Model reduction (gyroscope mass matrix) \\
49 & \texttt{bcsstk17} & 10\,974 & 428\,650 & Structural (elevated pressure vessel) \\
50 & \texttt{wathen100} & 30\,401 & 471\,601 & Structural (Wathen FE) \\
51 & \texttt{gridgena} & 48\,962 & 512\,084 & Optimisation (grid generation) \\
52 & \texttt{apache1} & 80\,800 & 542\,184 & 3D structural mechanics \\
53 & \texttt{Trefethen\_20000} & 20\,000 & 554\,466 & Combinatorics (Trefethen) \\
54 & \texttt{wathen120} & 36\,441 & 565\,761 & Structural (Wathen FE) \\
55 & \texttt{thermal1} & 82\,654 & 574\,458 & Steady-state thermal \\
56 & \texttt{crystm03} & 24\,696 & 583\,770 & Materials (crystal mass matrix) \\
57 & \texttt{finan512} & 74\,752 & 596\,992 & Financial portfolio optimisation \\
58 & \texttt{cbuckle} & 13\,681 & 676\,515 & Structural (buckling) \\
59 & \texttt{thermomech\_TC} & 102\,158 & 711\,558 & Thermomechanical \\
60 & \texttt{thermomech\_TK} & 102\,158 & 711\,558 & Thermomechanical \\
61 & \texttt{Pres\_Poisson} & 14\,822 & 715\,804 & CFD (pressure Poisson) \\
62 & \texttt{bundle1} & 10\,581 & 770\,901 & Computer vision (bundle adjustment) \\
63 & \texttt{bcsstk36} & 23\,052 & 1\,143\,140 & Structural (automobile shock absorber) \\
64 & \texttt{cfd1} & 70\,656 & 1\,828\,364 & CFD pressure matrix \\
65 & \texttt{bcsstk39} & 46\,772 & 2\,089\,294 & Structural (solid element model) \\
66 & \texttt{nd3k} & 9\,000 & 3\,279\,690 & 3D mesh problem (ND set) \\
67 & \texttt{parabolic\_fem} & 525\,825 & 3\,674\,625 & Parabolic FEM, convection-diffusion \\
68 & \texttt{offshore} & 259\,789 & 4\,242\,673 & Electromagnetics (offshore EM survey) \\
69 & \texttt{apache2} & 715\,176 & 4\,817\,870 & 3D structural mechanics \\
70 & \texttt{ecology2} & 999\,999 & 4\,995\,991 & 2D circular food-web model \\
71 & \texttt{G3\_circuit} & 1\,585\,478 & 7\,660\,826 & Circuit simulation \\
72 & \texttt{thermal2} & 1\,228\,045 & 8\,580\,313 & Steady-state thermal \\
73 & \texttt{crankseg\_1} & 52\,804 & 10\,614\,210 & Structural (crank segment) \\
74 & \texttt{crankseg\_2} & 63\,838 & 14\,148\,858 & Structural (crank segment) \\
75 & \texttt{af\_shell3} & 504\,855 & 17\,588\,875 & Structural (sheet-metal forming) \\
76 & \texttt{StocF-1465} & 1\,465\,137 & 21\,005\,389 & Fluid flow (stochastic FEM) \\
\end{longtable}

\begin{longtable}{rlrrl}
\caption{The nonsymmetric SuiteSparse/Matrix Market test matrices (numbers 77 to 91), with the size $n$ and number of nonzeros (nnz) of the full matrix and the application area recorded by the collections. The numbering is that used throughout the article.}
\label{tab:matrices_asym} \\
\toprule
\# & Matrix & $n$ & nnz & Application area \\
\midrule
\endfirsthead
\toprule
\# & Matrix & $n$ & nnz & Application area \\
\midrule
\endhead
\midrule
\endfoot
\bottomrule
\endlastfoot
77 & \texttt{cage10} & 11\,397 & 150\,645 & DNA electrophoresis \\
78 & \texttt{cage11} & 39\,082 & 559\,722 & DNA electrophoresis \\
79 & \texttt{ASIC\_100ks} & 99\,190 & 578\,890 & Circuit simulation \\
80 & \texttt{trans4} & 116\,835 & 766\,396 & Circuit simulation \\
81 & \texttt{torso2} & 115\,967 & 1\,033\,473 & Bioengineering (electrocardiology) \\
82 & \texttt{3D\_51448\_3D} & 51\,448 & 1\,056\,610 & Structural (3D problem) \\
83 & \texttt{venkat01} & 62\,424 & 1\,717\,792 & CFD (unstructured 2D Euler) \\
84 & \texttt{majorbasis} & 160\,000 & 1\,750\,416 & Optimisation \\
85 & \texttt{poisson3Db} & 85\,623 & 2\,374\,949 & CFD (3D Poisson) \\
86 & \texttt{stomach} & 213\,360 & 3\,021\,648 & Bioengineering (3D EM model of duodenum) \\
87 & \texttt{tmt\_unsym} & 917\,825 & 4\,584\,801 & Electromagnetics \\
88 & \texttt{atmosmodd} & 1\,270\,432 & 8\,814\,880 & Atmospheric modelling \\
89 & \texttt{Freescale1} & 3\,428\,755 & 18\,920\,347 & Circuit simulation \\
90 & \texttt{CoupCons3D} & 416\,800 & 22\,322\,336 & Structural (coupled consolidation) \\
91 & \texttt{Transport} & 1\,602\,111 & 23\,500\,731 & Structural (3D transport) \\
\end{longtable}
\endgroup

\section{Full results}
\label{sec:sm_results}
\Cref{tab:results_sym_ilu0,tab:results_sym_ilu1,tab:results_asym_ilu0,tab:results_asym_ilu1} give the per-matrix iteration counts for the $\mat{L}$ and $\mat{U}$ factor solves discussed in \cref{sec:lu_solves} of the article, for the symmetric and nonsymmetric matrices with ILU(0) and ILU(1) factorisations. The solver parameters are those of \cref{sec:parameters}.

The ordering column of each table gives the winning ordering selected by the sweep described in \cref{sec:ilu} of the article. The orderings swept are \textsc{natural}, \textsc{nd}, \textsc{rcm}, \textsc{qmd}, \textsc{1wd} and \textsc{rowlength} from PETSc, and \textsc{colamd}, \textsc{mmd\_ata} and \textsc{mmd\_at\_plus\_a} from SciPy. The last three are computed from the column structure (of $\mat{A}$, $\mat{A}^\textrm{T}\mat{A}$ and $\mat{A}^\textrm{T}+\mat{A}$, respectively), and SciPy exposes them only through its SuperLU interface, which returns a row permutation alongside the column ordering, chosen by its own threshold partial pivoting. We therefore extract these orderings from a probe factorisation with pivoting disabled (a diagonal pivot threshold of zero, so that the diagonal entry is always accepted as the pivot), in which case the row and column permutations returned coincide and the reordering is symmetric, $\mat{P}\mat{A}\mat{P}^\textrm{T}$.

\par\begingroup\footnotesize
\newcolumntype{Z}{>{\centering\arraybackslash}p{2.5em}}
\begin{longtable}{r l l Z@{\hskip 0.6em}Z!{\color{light}\vrule}Z@{\hskip 0.6em}Z!{\color{light}\vrule}Z@{\hskip 0.6em}Z!{\color{light}\vrule}Z@{\hskip 0.6em}Z}
\caption{For the symmetric matrices from the SuiteSparse collection and Matrix Market (\cref{tab:matrices_sym}) the number of iterations to solve the $\mat{L}$ and $\mat{U}$ factors generated by an ILU(0) factorisation to a relative tolerance of $10^{-6}$. ``MAX'' denotes hitting the maximum iteration count (2000), ``DIV'' indicates the solve diverged and ``DNF'' indicates the method did not finish, either running out of memory or exceeding the time limit.}
\label{tab:results_sym_ilu0} \\
\toprule
\# & Matrix & ordering & \multicolumn{2}{c}{AIRG} & \multicolumn{2}{!{\color{light}\vrule}c}{GMRES poly} & \multicolumn{2}{!{\color{light}\vrule}c}{ISAI} & \multicolumn{2}{!{\color{light}\vrule}c}{Jacobi} \\
\cmidrule(lr){4-5}\cmidrule(lr){6-7}\cmidrule(lr){8-9}\cmidrule(lr){10-11}
 & & & L & U & L & U & L & U & L & U \\
\midrule
\endfirsthead
\toprule
\# & Matrix & ordering & \multicolumn{2}{c}{AIRG} & \multicolumn{2}{!{\color{light}\vrule}c}{GMRES poly} & \multicolumn{2}{!{\color{light}\vrule}c}{ISAI} & \multicolumn{2}{!{\color{light}\vrule}c}{Jacobi} \\
\cmidrule(lr){4-5}\cmidrule(lr){6-7}\cmidrule(lr){8-9}\cmidrule(lr){10-11}
 & & & L & U & L & U & L & U & L & U \\
\midrule
\endhead
\midrule
\endfoot
\bottomrule
\endlastfoot
1 & \texttt{1138\_bus} & \texttt{mmd\_ata} & 2 & 2 & 2 & 2 & 5 & 5 & 9 & 10 \\
2 & \texttt{bcsstk19} & \texttt{natural} & 5 & 4 & \cellcolor{light}MAX & \cellcolor{light}MAX & 63 & 68 & 155 & 157 \\
3 & \texttt{Chem97ZtZ} & \texttt{natural} & 1 & 1 & 1 & 1 & 1 & 1 & 2 & 2 \\
4 & \texttt{bcsstk08} & \texttt{1wd} & 2 & 3 & 4 & 12 & 7 & 7 & 17 & 16 \\
5 & \texttt{mhd3200b} & \texttt{rcm} & 6 & 6 & 4 & 5 & 8 & 9 & 21 & 22 \\
6 & \texttt{bcsstm12} & \texttt{nd} & 3 & 4 & 11 & 18 & 5 & 5 & 33 & 34 \\
7 & \texttt{bcsstk10} & \texttt{rcm} & 4 & 4 & 94 & 40 & 29 & 29 & \cellcolor{light}DIV & \cellcolor{light}DIV \\
8 & \texttt{mhd4800b} & \texttt{rcm} & 6 & 6 & 4 & 7 & 9 & 9 & 22 & 23 \\
9 & \texttt{msc01050} & \texttt{natural} & 3 & 3 & 761 & 1219 & 8 & 8 & 24 & 24 \\
10 & \texttt{plbuckle} & \texttt{natural} & 3 & 5 & 16 & 78 & 15 & 15 & 41 & 42 \\
11 & \texttt{Trefethen\_2000} & \texttt{natural} & 2 & 2 & 1 & 2 & 3 & 4 & 6 & 8 \\
12 & \texttt{bcsstk27} & \texttt{natural} & 3 & 3 & 22 & 19 & 20 & 20 & \cellcolor{light}DIV & \cellcolor{light}DIV \\
13 & \texttt{bcsstk14} & \texttt{rowlength} & 2 & 3 & 3 & 3 & 7 & 7 & 16 & 16 \\
14 & \texttt{nasa2146} & \texttt{natural} & 3 & 3 & 17 & 14 & 20 & 20 & 62 & 60 \\
15 & \texttt{fv1} & \texttt{natural} & 2 & 2 & 3 & 3 & 10 & 10 & 22 & 22 \\
16 & \texttt{fv2} & \texttt{natural} & 2 & 2 & 3 & 3 & 10 & 10 & 22 & 22 \\
17 & \texttt{msc04515} & \texttt{nd} & 3 & 5 & \cellcolor{light}MAX & \cellcolor{light}MAX & 20 & 18 & \cellcolor{light}DIV & \cellcolor{light}DIV \\
18 & \texttt{crystm01} & \texttt{natural} & 3 & 3 & 5 & 5 & 9 & 9 & 36 & 36 \\
19 & \texttt{bcsstk15} & \texttt{colamd} & 3 & 3 & 7 & 7 & 8 & 9 & 20 & 19 \\
20 & \texttt{bodyy4} & \texttt{natural} & 3 & 2 & 9 & 7 & 13 & 15 & 26 & 30 \\
21 & \texttt{aft01} & \texttt{rcm} & 4 & 4 & 23 & 25 & 108 & \cellcolor{light}DIV & \cellcolor{light}DIV & \cellcolor{light}DIV \\
22 & \texttt{bodyy5} & \texttt{natural} & 2 & 2 & 20 & 11 & 15 & 18 & 30 & 37 \\
23 & \texttt{bodyy6} & \texttt{natural} & 3 & 2 & 33 & 14 & 19 & 21 & 38 & 41 \\
24 & \texttt{ted\_B} & \texttt{natural} & 3 & 2 & 3 & 3 & 5 & 5 & 17 & 18 \\
25 & \texttt{ted\_B\_unscaled} & \texttt{natural} & 3 & 2 & 3 & 3 & 5 & 5 & 17 & 18 \\
26 & \texttt{bcsstk18} & \texttt{natural} & 3 & 3 & 6 & 20 & 11 & 12 & 25 & 24 \\
27 & \texttt{bcsstk24} & \texttt{rcm} & 5 & 5 & \cellcolor{light}MAX & \cellcolor{light}MAX & 19 & 17 & \cellcolor{light}DIV & 107 \\
28 & \texttt{Muu} & \texttt{rcm} & 3 & 3 & 3 & 3 & 6 & 6 & 19 & 19 \\
29 & \texttt{t2dah\_e} & \texttt{rcm} & 4 & 5 & 5 & 5 & 52 & 58 & \cellcolor{light}DIV & \cellcolor{light}DIV \\
30 & \texttt{obstclae} & \texttt{natural} & 2 & 2 & 4 & 4 & 12 & 12 & 23 & 23 \\
31 & \texttt{torsion1} & \texttt{natural} & 2 & 2 & 4 & 4 & 12 & 12 & 23 & 23 \\
32 & \texttt{jnlbrng1} & \texttt{natural} & 2 & 2 & 7 & 7 & 19 & 19 & 37 & 37 \\
33 & \texttt{minsurfo} & \texttt{natural} & 3 & 3 & 4 & 4 & 12 & 12 & 24 & 24 \\
34 & \texttt{s3rmt3m3} & \texttt{rcm} & 9 & 7 & \cellcolor{light}MAX & \cellcolor{light}DIV & 33 & 30 & \cellcolor{light}DIV & \cellcolor{light}DIV \\
35 & \texttt{s1rmt3m1} & \texttt{natural} & 4 & 4 & 17 & 24 & 18 & 18 & 74 & 74 \\
36 & \texttt{s2rmt3m1} & \texttt{natural} & 5 & 4 & 1619 & 51 & 19 & 18 & 87 & 86 \\
37 & \texttt{s3rmt3m1} & \texttt{natural} & 6 & 9 & \cellcolor{light}MAX & \cellcolor{light}MAX & 17 & 17 & \cellcolor{light}DIV & \cellcolor{light}DIV \\
38 & \texttt{Dubcova1} & \texttt{rcm} & 3 & 3 & 7 & 7 & 16 & 16 & 49 & 49 \\
39 & \texttt{s1rmq4m1} & \texttt{natural} & 4 & 4 & 35 & 64 & 20 & 21 & 83 & 86 \\
40 & \texttt{s2rmq4m1} & \texttt{natural} & 6 & 5 & \cellcolor{light}MAX & 77 & 20 & 20 & 75 & 71 \\
41 & \texttt{s3rmq4m1} & \texttt{natural} & 9 & 9 & \cellcolor{light}DIV & \cellcolor{light}DIV & 20 & 20 & \cellcolor{light}DIV & \cellcolor{light}DIV \\
42 & \texttt{chipcool0} & \texttt{1wd} & 3 & 3 & 4 & 4 & 13 & 12 & 29 & 28 \\
43 & \texttt{bcsstk16} & \texttt{natural} & 3 & 3 & 6 & 7 & 18 & 19 & 41 & 39 \\
44 & \texttt{crystm02} & \texttt{natural} & 3 & 3 & 5 & 5 & 9 & 9 & 47 & 47 \\
45 & \texttt{shallow\_water1} & \texttt{natural} & 2 & 2 & 2 & 2 & 5 & 5 & 9 & 9 \\
46 & \texttt{shallow\_water2} & \texttt{rcm} & 2 & 2 & 3 & 3 & 9 & 9 & 17 & 17 \\
47 & \texttt{Kuu} & \texttt{rowlength} & 3 & 3 & 3 & 3 & 6 & 5 & 16 & 15 \\
48 & \texttt{gyro\_m} & \texttt{1wd} & 4 & 3 & 7 & 7 & 33 & 32 & \cellcolor{light}DIV & \cellcolor{light}DIV \\
49 & \texttt{bcsstk17} & \texttt{nd} & 3 & 3 & 109 & 86 & 10 & 9 & 59 & \cellcolor{light}DIV \\
50 & \texttt{wathen100} & \texttt{natural} & 4 & 4 & 4 & 4 & 44 & 45 & \cellcolor{light}DIV & \cellcolor{light}DIV \\
51 & \texttt{gridgena} & \texttt{rcm} & 3 & 3 & 6 & 6 & 15 & 15 & 42 & 39 \\
52 & \texttt{apache1} & \texttt{rcm} & 2 & 3 & 3 & 4 & 9 & 10 & 18 & 20 \\
53 & \texttt{Trefethen\_20000} & \texttt{colamd} & 2 & 1 & 1 & 1 & 3 & 3 & 6 & 6 \\
54 & \texttt{wathen120} & \texttt{natural} & 4 & 4 & 4 & 4 & 45 & 45 & \cellcolor{light}DIV & \cellcolor{light}DIV \\
55 & \texttt{thermal1} & \texttt{rcm} & 3 & 3 & 5 & 5 & 14 & 14 & 31 & 31 \\
56 & \texttt{crystm03} & \texttt{natural} & 3 & 3 & 5 & 5 & 9 & 9 & 52 & 52 \\
57 & \texttt{finan512} & \texttt{natural} & 2 & 2 & 2 & 2 & 4 & 4 & 8 & 8 \\
58 & \texttt{cbuckle} & \texttt{qmd} & 3 & 4 & 77 & 141 & 8 & 8 & \cellcolor{light}DIV & \cellcolor{light}DIV \\
59 & \texttt{thermomech\_TC} & \texttt{rcm} & 2 & 2 & 2 & 2 & 8 & 8 & 18 & 19 \\
60 & \texttt{thermomech\_TK} & \texttt{rcm} & 3 & 3 & 5 & 5 & 15 & 15 & 32 & 32 \\
61 & \texttt{Pres\_Poisson} & \texttt{rcm} & 4 & 4 & 18 & 11 & 20 & 19 & \cellcolor{light}DIV & \cellcolor{light}DIV \\
62 & \texttt{bundle1} & \texttt{qmd} & 2 & 2 & 2 & 2 & 3 & 3 & 10 & 9 \\
63 & \texttt{bcsstk36} & \texttt{1wd} & \cellcolor{light}DIV & \cellcolor{light}DIV & \cellcolor{light}DIV & \cellcolor{light}MAX & \cellcolor{light}DIV & \cellcolor{light}DIV & \cellcolor{light}DIV & \cellcolor{light}DIV \\
64 & \texttt{cfd1} & \texttt{1wd} & 5 & 5 & 31 & 31 & 56 & 56 & \cellcolor{light}DIV & \cellcolor{light}DIV \\
65 & \texttt{bcsstk39} & \texttt{mmd\_ata} & 5 & 6 & \cellcolor{light}MAX & \cellcolor{light}MAX & 31 & 31 & \cellcolor{light}DIV & \cellcolor{light}DIV \\
66 & \texttt{nd3k} & \texttt{natural} & 4 & 4 & 265 & 217 & \cellcolor{light}DIV & \cellcolor{light}DIV & \cellcolor{light}DIV & \cellcolor{light}DIV \\
67 & \texttt{parabolic\_fem} & \texttt{rcm} & 3 & 3 & 8 & 8 & 22 & 22 & 46 & 46 \\
68 & \texttt{offshore} & \texttt{1wd} & 3 & 3 & 5 & 4 & 13 & 11 & 30 & 26 \\
69 & \texttt{apache2} & \texttt{natural} & 2 & 2 & 4 & 4 & 12 & 12 & 24 & 24 \\
70 & \texttt{ecology2} & \texttt{rcm} & 3 & 3 & 5 & 5 & 13 & 13 & 26 & 26 \\
71 & \texttt{G3\_circuit} & \texttt{natural} & 3 & 3 & 10 & 9 & 31 & 29 & 62 & 57 \\
72 & \texttt{thermal2} & \texttt{1wd} & 3 & 3 & 4 & 4 & 13 & 13 & 29 & 29 \\
73 & \texttt{crankseg\_1} & \texttt{rcm} & 3 & 3 & 21 & 19 & 22 & 19 & \cellcolor{light}DIV & \cellcolor{light}DIV \\
74 & \texttt{crankseg\_2} & \texttt{rcm} & 3 & \cellcolor{light}DNF & 24 & 22 & 25 & 22 & \cellcolor{light}DIV & \cellcolor{light}DIV \\
75 & \texttt{af\_shell3} & \texttt{rcm} & 3 & 3 & 8 & 8 & 26 & 26 & 68 & 68 \\
76 & \texttt{StocF-1465} & \texttt{rcm} & 4 & 2 & \cellcolor{light}MAX & 21 & 35 & 33 & 73 & 68 \\
\end{longtable}

\begin{longtable}{r l l Z@{\hskip 0.6em}Z!{\color{light}\vrule}Z@{\hskip 0.6em}Z!{\color{light}\vrule}Z@{\hskip 0.6em}Z!{\color{light}\vrule}Z@{\hskip 0.6em}Z}
\caption{For the symmetric matrices from the SuiteSparse collection and Matrix Market (\cref{tab:matrices_sym}) the number of iterations to solve the $\mat{L}$ and $\mat{U}$ factors generated by an ILU(1) factorisation to a relative tolerance of $10^{-6}$. ``MAX'' denotes hitting the maximum iteration count (2000), ``DIV'' indicates the solve diverged and ``DNF'' indicates the method did not finish, either running out of memory or exceeding the time limit.}
\label{tab:results_sym_ilu1} \\
\toprule
\# & Matrix & ordering & \multicolumn{2}{c}{AIRG} & \multicolumn{2}{!{\color{light}\vrule}c}{GMRES poly} & \multicolumn{2}{!{\color{light}\vrule}c}{ISAI} & \multicolumn{2}{!{\color{light}\vrule}c}{Jacobi} \\
\cmidrule(lr){4-5}\cmidrule(lr){6-7}\cmidrule(lr){8-9}\cmidrule(lr){10-11}
 & & & L & U & L & U & L & U & L & U \\
\midrule
\endfirsthead
\toprule
\# & Matrix & ordering & \multicolumn{2}{c}{AIRG} & \multicolumn{2}{!{\color{light}\vrule}c}{GMRES poly} & \multicolumn{2}{!{\color{light}\vrule}c}{ISAI} & \multicolumn{2}{!{\color{light}\vrule}c}{Jacobi} \\
\cmidrule(lr){4-5}\cmidrule(lr){6-7}\cmidrule(lr){8-9}\cmidrule(lr){10-11}
 & & & L & U & L & U & L & U & L & U \\
\midrule
\endhead
\midrule
\endfoot
\bottomrule
\endlastfoot
1 & \texttt{1138\_bus} & \texttt{rcm} & 3 & 2 & 238 & 7 & 17 & 16 & 37 & 32 \\
2 & \texttt{bcsstk19} & \texttt{rcm} & 4 & 4 & \cellcolor{light}MAX & \cellcolor{light}MAX & 45 & 45 & \cellcolor{light}DIV & \cellcolor{light}DIV \\
3 & \texttt{Chem97ZtZ} & \texttt{natural} & 1 & 1 & 1 & 1 & 1 & 1 & 2 & 2 \\
4 & \texttt{bcsstk08} & \texttt{rcm} & 3 & 3 & 5 & 21 & 7 & 7 & 16 & 17 \\
5 & \texttt{mhd3200b} & \texttt{natural} & 6 & 6 & 5 & 6 & 8 & 9 & 19 & 22 \\
6 & \texttt{bcsstm12} & \texttt{nd} & 2 & 3 & 11 & 17 & 4 & 4 & 33 & 34 \\
7 & \texttt{bcsstk10} & \texttt{1wd} & 3 & 3 & 12 & 26 & 6 & 6 & 28 & 28 \\
8 & \texttt{mhd4800b} & \texttt{natural} & 6 & 6 & 5 & 7 & 9 & 9 & 20 & 22 \\
9 & \texttt{msc01050} & \texttt{qmd} & 2 & 3 & 8 & 10 & 7 & 8 & 31 & 32 \\
10 & \texttt{plbuckle} & \texttt{natural} & 4 & 5 & 10 & 21 & 16 & 16 & 43 & 46 \\
11 & \texttt{Trefethen\_2000} & \texttt{nd} & 2 & 2 & 1 & 1 & 3 & 3 & 6 & 6 \\
12 & \texttt{bcsstk27} & \texttt{natural} & 4 & 4 & 36 & 23 & 21 & 21 & \cellcolor{light}DIV & \cellcolor{light}DIV \\
13 & \texttt{bcsstk14} & \texttt{rcm} & 3 & 3 & 11 & 21 & 16 & 16 & 45 & 50 \\
14 & \texttt{nasa2146} & \texttt{natural} & 3 & 3 & 18 & 14 & 17 & 17 & 61 & 59 \\
15 & \texttt{fv1} & \texttt{natural} & 2 & 2 & 4 & 4 & 11 & 11 & 26 & 26 \\
16 & \texttt{fv2} & \texttt{natural} & 2 & 2 & 4 & 4 & 11 & 11 & 26 & 26 \\
17 & \texttt{msc04515} & \texttt{rcm} & 5 & 4 & 138 & 46 & 28 & 26 & \cellcolor{light}DIV & \cellcolor{light}DIV \\
18 & \texttt{crystm01} & \texttt{natural} & 3 & 3 & 5 & 5 & 6 & 6 & 36 & 36 \\
19 & \texttt{bcsstk15} & \texttt{rcm} & 4 & 4 & 44 & 20 & 14 & 14 & 36 & 38 \\
20 & \texttt{bodyy4} & \texttt{natural} & 3 & 3 & 9 & 7 & 13 & 15 & 26 & 30 \\
21 & \texttt{aft01} & \texttt{natural} & 6 & 3 & 29 & 30 & 78 & 79 & \cellcolor{light}DIV & \cellcolor{light}DIV \\
22 & \texttt{bodyy5} & \texttt{natural} & 3 & 2 & 20 & 11 & 15 & 18 & 30 & 36 \\
23 & \texttt{bodyy6} & \texttt{natural} & 3 & 2 & 33 & 14 & 19 & 21 & 38 & 41 \\
24 & \texttt{ted\_B} & \texttt{natural} & 3 & 2 & 3 & 3 & 5 & 5 & 17 & 18 \\
25 & \texttt{ted\_B\_unscaled} & \texttt{natural} & 3 & 2 & 3 & 3 & 5 & 5 & 17 & 18 \\
26 & \texttt{bcsstk18} & \texttt{mmd\_at\_plus\_a} & 3 & 3 & 95 & 7 & 8 & 7 & 19 & 17 \\
27 & \texttt{bcsstk24} & \texttt{rcm} & 6 & 5 & \cellcolor{light}MAX & \cellcolor{light}MAX & 27 & 26 & \cellcolor{light}DIV & \cellcolor{light}DIV \\
28 & \texttt{Muu} & \texttt{rcm} & 3 & 3 & 3 & 3 & 6 & 6 & 19 & 19 \\
29 & \texttt{t2dah\_e} & \texttt{natural} & 3 & 4 & 3 & 3 & 7 & 7 & 95 & 144 \\
30 & \texttt{obstclae} & \texttt{natural} & 3 & 3 & 6 & 6 & 17 & 17 & 35 & 35 \\
31 & \texttt{torsion1} & \texttt{natural} & 3 & 3 & 6 & 6 & 17 & 17 & 35 & 35 \\
32 & \texttt{jnlbrng1} & \texttt{natural} & 3 & 3 & 12 & 12 & 33 & 33 & 66 & 66 \\
33 & \texttt{minsurfo} & \texttt{natural} & 3 & 3 & 6 & 6 & 18 & 18 & 37 & 37 \\
34 & \texttt{s3rmt3m3} & \texttt{rcm} & 9 & 8 & \cellcolor{light}MAX & \cellcolor{light}DIV & 36 & 32 & \cellcolor{light}DIV & \cellcolor{light}DIV \\
35 & \texttt{s1rmt3m1} & \texttt{rcm} & 4 & 4 & 80 & 30 & 32 & 31 & 91 & 90 \\
36 & \texttt{s2rmt3m1} & \texttt{natural} & 4 & 4 & \cellcolor{light}MAX & 171 & 23 & 21 & 93 & 92 \\
37 & \texttt{s3rmt3m1} & \texttt{natural} & 8 & 10 & \cellcolor{light}MAX & \cellcolor{light}DIV & 22 & 20 & \cellcolor{light}DIV & \cellcolor{light}DIV \\
38 & \texttt{Dubcova1} & \texttt{rcm} & 3 & 3 & 9 & 9 & 22 & 20 & 51 & 50 \\
39 & \texttt{s1rmq4m1} & \texttt{natural} & 4 & 4 & 84 & 89 & 22 & 22 & 82 & 83 \\
40 & \texttt{s2rmq4m1} & \texttt{natural} & 6 & 6 & \cellcolor{light}MAX & 170 & 23 & 22 & 106 & 102 \\
41 & \texttt{s3rmq4m1} & \texttt{natural} & 9 & 8 & \cellcolor{light}MAX & \cellcolor{light}DIV & 24 & 23 & \cellcolor{light}DIV & \cellcolor{light}DIV \\
42 & \texttt{chipcool0} & \texttt{natural} & 3 & 3 & 7 & 7 & 20 & 20 & 50 & 50 \\
43 & \texttt{bcsstk16} & \texttt{natural} & 3 & 3 & 9 & 9 & 22 & 23 & 53 & 53 \\
44 & \texttt{crystm02} & \texttt{natural} & 3 & 3 & 5 & 5 & 6 & 6 & 47 & 47 \\
45 & \texttt{shallow\_water1} & \texttt{rcm} & 2 & 2 & 2 & 2 & 6 & 6 & 11 & 11 \\
46 & \texttt{shallow\_water2} & \texttt{rcm} & 2 & 3 & 3 & 3 & 10 & 10 & 19 & 19 \\
47 & \texttt{Kuu} & \texttt{1wd} & 3 & 3 & 17 & 17 & 22 & 23 & 96 & 92 \\
48 & \texttt{gyro\_m} & \texttt{natural} & 3 & 3 & 13 & 11 & 21 & 20 & \cellcolor{light}DIV & \cellcolor{light}DIV \\
49 & \texttt{bcsstk17} & \texttt{rowlength} & 3 & 4 & 365 & 76 & 15 & 14 & \cellcolor{light}DIV & \cellcolor{light}DIV \\
50 & \texttt{wathen100} & \texttt{natural} & 4 & 4 & 4 & 4 & 21 & 21 & \cellcolor{light}DIV & \cellcolor{light}DIV \\
51 & \texttt{gridgena} & \texttt{qmd} & 2 & 2 & 3 & 3 & 8 & 9 & 19 & 21 \\
52 & \texttt{apache1} & \texttt{rcm} & 3 & 3 & 4 & 5 & 13 & 14 & 27 & 30 \\
53 & \texttt{Trefethen\_20000} & \texttt{natural} & 2 & 2 & 1 & 2 & 3 & 4 & 6 & 8 \\
54 & \texttt{wathen120} & \texttt{natural} & 4 & 4 & 4 & 4 & 21 & 21 & \cellcolor{light}DIV & \cellcolor{light}DIV \\
55 & \texttt{thermal1} & \texttt{1wd} & 3 & 3 & 7 & 6 & 19 & 19 & 44 & 44 \\
56 & \texttt{crystm03} & \texttt{natural} & 3 & 3 & 5 & 5 & 6 & 6 & 52 & 52 \\
57 & \texttt{finan512} & \texttt{natural} & 2 & 2 & 2 & 2 & 4 & 4 & 8 & 8 \\
58 & \texttt{cbuckle} & \texttt{rcm} & 4 & 4 & 794 & \cellcolor{light}MAX & 22 & 23 & \cellcolor{light}DIV & \cellcolor{light}DIV \\
59 & \texttt{thermomech\_TC} & \texttt{rcm} & 2 & 2 & 2 & 2 & 7 & 7 & 17 & 17 \\
60 & \texttt{thermomech\_TK} & \texttt{rcm} & 3 & 3 & 7 & 7 & 22 & 22 & 50 & 49 \\
61 & \texttt{Pres\_Poisson} & \texttt{rcm} & 5 & 3 & 20 & 14 & 26 & 25 & \cellcolor{light}DIV & \cellcolor{light}DIV \\
62 & \texttt{bundle1} & \texttt{rowlength} & 2 & 2 & 2 & 2 & 3 & 3 & 12 & 12 \\
63 & \texttt{bcsstk36} & \texttt{rcm} & 8 & 9 & \cellcolor{light}DIV & \cellcolor{light}DIV & \cellcolor{light}DIV & \cellcolor{light}DIV & \cellcolor{light}DIV & \cellcolor{light}DIV \\
64 & \texttt{cfd1} & \texttt{1wd} & 5 & 5 & 49 & 51 & 53 & 54 & \cellcolor{light}DIV & \cellcolor{light}DIV \\
65 & \texttt{bcsstk39} & \texttt{nd} & 4 & 5 & \cellcolor{light}MAX & \cellcolor{light}MAX & \cellcolor{light}DIV & \cellcolor{light}DIV & \cellcolor{light}DIV & \cellcolor{light}DIV \\
66 & \texttt{nd3k} & \texttt{natural} & 3 & 3 & 133 & 95 & \cellcolor{light}DIV & \cellcolor{light}DIV & \cellcolor{light}DIV & \cellcolor{light}DIV \\
67 & \texttt{parabolic\_fem} & \texttt{rcm} & 4 & 3 & 11 & 11 & 31 & 31 & 66 & 66 \\
68 & \texttt{offshore} & \texttt{natural} & 4 & 4 & 40 & 64 & 28 & 29 & 83 & \cellcolor{light}DIV \\
69 & \texttt{apache2} & \texttt{natural} & 3 & 3 & 6 & 6 & 18 & 19 & 38 & 40 \\
70 & \texttt{ecology2} & \texttt{natural} & 3 & 3 & 8 & 8 & 22 & 22 & 46 & 46 \\
71 & \texttt{G3\_circuit} & \texttt{natural} & 4 & 3 & 10 & 9 & 32 & 30 & 63 & 59 \\
72 & \texttt{thermal2} & \texttt{rcm} & 3 & 3 & 7 & 7 & 21 & 21 & 49 & 48 \\
73 & \texttt{crankseg\_1} & \texttt{rcm} & \cellcolor{light}DNF & \cellcolor{light}DNF & 22 & 20 & 19 & 18 & \cellcolor{light}DIV & \cellcolor{light}DIV \\
74 & \texttt{crankseg\_2} & \texttt{rcm} & \cellcolor{light}DNF & \cellcolor{light}DNF & 26 & 23 & 21 & 19 & \cellcolor{light}DIV & \cellcolor{light}DIV \\
75 & \texttt{af\_shell3} & \texttt{rcm} & 3 & 3 & 10 & 10 & 32 & 32 & 96 & 94 \\
76 & \texttt{StocF-1465} & \texttt{rcm} & 5 & 3 & \cellcolor{light}MAX & 23 & 41 & 37 & 94 & 80 \\
\end{longtable}

\begin{longtable}{r l l Z@{\hskip 0.6em}Z!{\color{light}\vrule}Z@{\hskip 0.6em}Z!{\color{light}\vrule}Z@{\hskip 0.6em}Z!{\color{light}\vrule}Z@{\hskip 0.6em}Z}
\caption{For the nonsymmetric matrices from the SuiteSparse collection and Matrix Market (\cref{tab:matrices_asym}) the number of iterations to solve the $\mat{L}$ and $\mat{U}$ factors generated by an ILU(0) factorisation to a relative tolerance of $10^{-6}$. ``MAX'' denotes hitting the maximum iteration count (2000), ``DIV'' indicates the solve diverged and ``DNF'' indicates the method did not finish, either running out of memory or exceeding the time limit.}
\label{tab:results_asym_ilu0} \\
\toprule
\# & Matrix & ordering & \multicolumn{2}{c}{AIRG} & \multicolumn{2}{!{\color{light}\vrule}c}{GMRES poly} & \multicolumn{2}{!{\color{light}\vrule}c}{ISAI} & \multicolumn{2}{!{\color{light}\vrule}c}{Jacobi} \\
\cmidrule(lr){4-5}\cmidrule(lr){6-7}\cmidrule(lr){8-9}\cmidrule(lr){10-11}
 & & & L & U & L & U & L & U & L & U \\
\midrule
\endfirsthead
\toprule
\# & Matrix & ordering & \multicolumn{2}{c}{AIRG} & \multicolumn{2}{!{\color{light}\vrule}c}{GMRES poly} & \multicolumn{2}{!{\color{light}\vrule}c}{ISAI} & \multicolumn{2}{!{\color{light}\vrule}c}{Jacobi} \\
\cmidrule(lr){4-5}\cmidrule(lr){6-7}\cmidrule(lr){8-9}\cmidrule(lr){10-11}
 & & & L & U & L & U & L & U & L & U \\
\midrule
\endhead
\midrule
\endfoot
\bottomrule
\endlastfoot
77 & \texttt{cage10} & \texttt{natural} & 2 & 2 & 2 & 2 & 5 & 5 & 10 & 11 \\
78 & \texttt{cage11} & \texttt{natural} & 2 & 2 & 2 & 2 & 5 & 5 & 10 & 11 \\
79 & \texttt{ASIC\_100ks} & \texttt{rcm} & 2 & 2 & 150 & 2 & 9 & 5 & 18 & 11 \\
80 & \texttt{trans4} & \texttt{rcm} & 2 & 2 & 1 & 2 & 3 & 4 & 7 & 8 \\
81 & \texttt{torso2} & \texttt{natural} & 2 & 2 & 3 & 3 & 7 & 7 & 16 & 16 \\
82 & \texttt{3D\_51448\_3D} & \texttt{rowlength} & 5 & 3 & \cellcolor{light}MAX & 6 & 16 & 16 & 57 & 32 \\
83 & \texttt{venkat01} & \texttt{rcm} & 5 & 4 & 11 & 14 & 29 & 33 & 59 & 82 \\
84 & \texttt{majorbasis} & \texttt{rcm} & 3 & 3 & 6 & 5 & 16 & 13 & 39 & 27 \\
85 & \texttt{poisson3Db} & \texttt{1wd} & 2 & 3 & 4 & 3 & 9 & 9 & 23 & 21 \\
86 & \texttt{stomach} & \texttt{rcm} & 3 & 3 & 5 & 5 & 14 & 14 & 27 & 28 \\
87 & \texttt{tmt\_unsym} & \texttt{rcm} & 3 & 3 & 5 & 5 & 13 & 13 & 26 & 26 \\
88 & \texttt{atmosmodd} & \texttt{rcm} & 2 & 2 & 4 & 4 & 13 & 12 & 25 & 24 \\
89 & \texttt{Freescale1} & \texttt{natural} & 2 & 2 & 11 & 11 & 33 & 32 & 66 & 63 \\
90 & \texttt{CoupCons3D} & \texttt{natural} & 3 & 4 & 6 & 7 & 17 & 18 & 37 & 41 \\
91 & \texttt{Transport} & \texttt{rcm} & 3 & 3 & 7 & 7 & 21 & 21 & 44 & 44 \\
\end{longtable}

\begin{longtable}{r l l Z@{\hskip 0.6em}Z!{\color{light}\vrule}Z@{\hskip 0.6em}Z!{\color{light}\vrule}Z@{\hskip 0.6em}Z!{\color{light}\vrule}Z@{\hskip 0.6em}Z}
\caption{For the nonsymmetric matrices from the SuiteSparse collection and Matrix Market (\cref{tab:matrices_asym}) the number of iterations to solve the $\mat{L}$ and $\mat{U}$ factors generated by an ILU(1) factorisation to a relative tolerance of $10^{-6}$. ``MAX'' denotes hitting the maximum iteration count (2000), ``DIV'' indicates the solve diverged and ``DNF'' indicates the method did not finish, either running out of memory or exceeding the time limit.}
\label{tab:results_asym_ilu1} \\
\toprule
\# & Matrix & ordering & \multicolumn{2}{c}{AIRG} & \multicolumn{2}{!{\color{light}\vrule}c}{GMRES poly} & \multicolumn{2}{!{\color{light}\vrule}c}{ISAI} & \multicolumn{2}{!{\color{light}\vrule}c}{Jacobi} \\
\cmidrule(lr){4-5}\cmidrule(lr){6-7}\cmidrule(lr){8-9}\cmidrule(lr){10-11}
 & & & L & U & L & U & L & U & L & U \\
\midrule
\endfirsthead
\toprule
\# & Matrix & ordering & \multicolumn{2}{c}{AIRG} & \multicolumn{2}{!{\color{light}\vrule}c}{GMRES poly} & \multicolumn{2}{!{\color{light}\vrule}c}{ISAI} & \multicolumn{2}{!{\color{light}\vrule}c}{Jacobi} \\
\cmidrule(lr){4-5}\cmidrule(lr){6-7}\cmidrule(lr){8-9}\cmidrule(lr){10-11}
 & & & L & U & L & U & L & U & L & U \\
\midrule
\endhead
\midrule
\endfoot
\bottomrule
\endlastfoot
77 & \texttt{cage10} & \texttt{natural} & 2 & 2 & 2 & 2 & 4 & 5 & 9 & 9 \\
78 & \texttt{cage11} & \texttt{natural} & 2 & 2 & 2 & 2 & 4 & 5 & 9 & 9 \\
79 & \texttt{ASIC\_100ks} & \texttt{natural} & 2 & 2 & 196 & 2 & 8 & 5 & 15 & 10 \\
80 & \texttt{trans4} & \texttt{rcm} & 2 & 2 & 2 & 2 & 4 & 4 & 13 & 13 \\
81 & \texttt{torso2} & \texttt{rcm} & 2 & 2 & 3 & 3 & 7 & 7 & 17 & 17 \\
82 & \texttt{3D\_51448\_3D} & \texttt{natural} & 5 & 4 & \cellcolor{light}MAX & 10 & 24 & 21 & \cellcolor{light}DIV & 44 \\
83 & \texttt{venkat01} & \texttt{rcm} & 5 & 4 & 13 & 16 & 34 & 36 & 70 & 89 \\
84 & \texttt{majorbasis} & \texttt{rcm} & 3 & 3 & 7 & 5 & 16 & 15 & 42 & 30 \\
85 & \texttt{poisson3Db} & \texttt{natural} & 2 & 2 & 3 & 3 & 6 & 6 & 19 & 19 \\
86 & \texttt{stomach} & \texttt{rcm} & 3 & 3 & 5 & 5 & 13 & 14 & 30 & 29 \\
87 & \texttt{tmt\_unsym} & \texttt{rcm} & 3 & 3 & 8 & 8 & 22 & 22 & 46 & 46 \\
88 & \texttt{atmosmodd} & \texttt{rcm} & 3 & 3 & 7 & 7 & 21 & 20 & 42 & 42 \\
89 & \texttt{Freescale1} & \texttt{natural} & 2 & 2 & 14 & 14 & 43 & 41 & 86 & 81 \\
90 & \texttt{CoupCons3D} & \texttt{natural} & 3 & 4 & 7 & 8 & 18 & 20 & 35 & 38 \\
91 & \texttt{Transport} & \texttt{rcm} & 3 & 3 & 11 & 11 & 30 & 30 & 64 & 64 \\
\end{longtable}
\endgroup

\section{Setup and solve time figures}
\label{sec:sm_figures}
\begin{figure}[ht]
\centering
\includegraphics[width=0.98\textwidth]{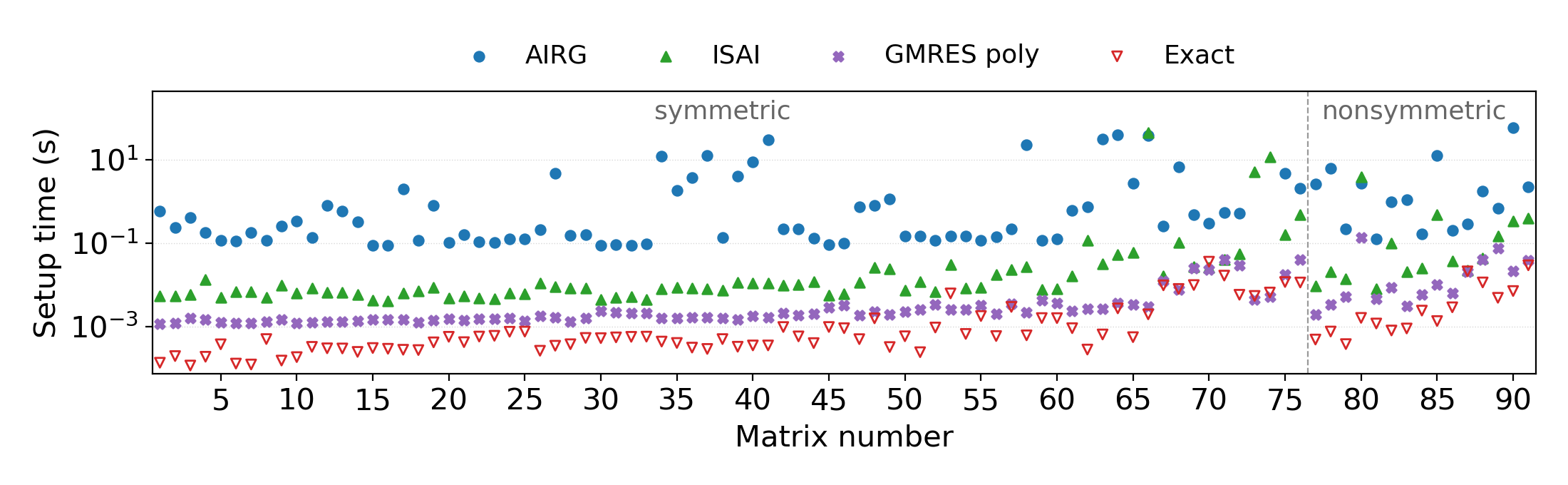}
\caption{Setup time comparison (smaller is better) between the GMRES polynomial, ISAI, AIRG and the analysis phase of the hipSPARSE exact triangular solve for the $\mat{L}$ and $\mat{U}$ factors generated by the ILU(1) factorisation of each test matrix, ordered as in the results tables (symmetric and nonsymmetric matrices separated by the dashed line). The results shown are taken as the maximum over the $\mat{L}$ and $\mat{U}$ factors. Jacobi requires no setup and is not shown. The ILU(0) factors are shown in \cref{fig:setup} of the article.}
\label{fig:setup_ilu1}
\end{figure}

\begin{figure}[ht]
\centering
\includegraphics[width=0.98\textwidth]{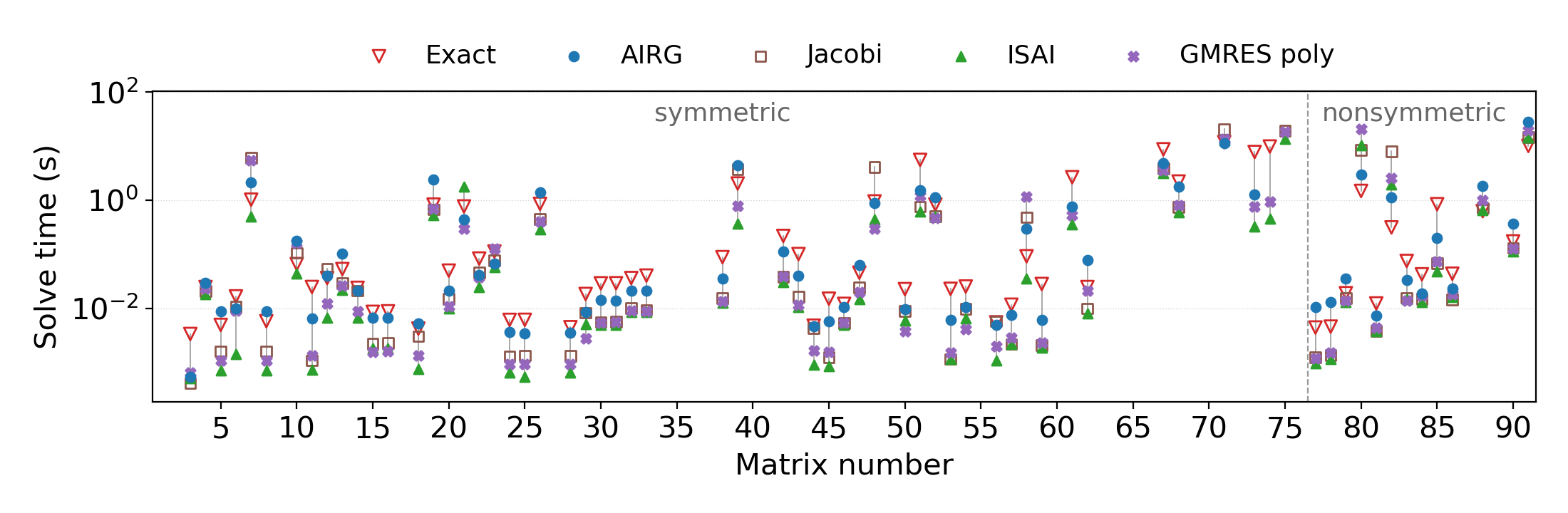}
\caption{Solve time comparison (smaller is better) for the $\mat{A}\mat{x}=\mat{b}$ solve with the ILU(0) factorisation for each test matrix, ordered as in the results tables (symmetric and nonsymmetric matrices separated by the dashed line). The solve uses an outer GMRES(30) preconditioned by the ILU factorisation, with the inner application of $\mat{L}^{-1}$ and $\mat{U}^{-1}$ performed by an exact triangular solve, by one V-cycle of AIRG, or by the Jacobi method, ISAI or the GMRES polynomial, each with the number of inner iterations that gave the smallest converged solve time for that problem. Only solves that converged to a relative tolerance of $10^{-6}$ are shown. \Cref{tab:axb_summary} of the article summarises these results.}
\label{fig:axb_timing_ilu0}
\end{figure}

\begin{figure}[ht]
\centering
\includegraphics[width=0.98\textwidth]{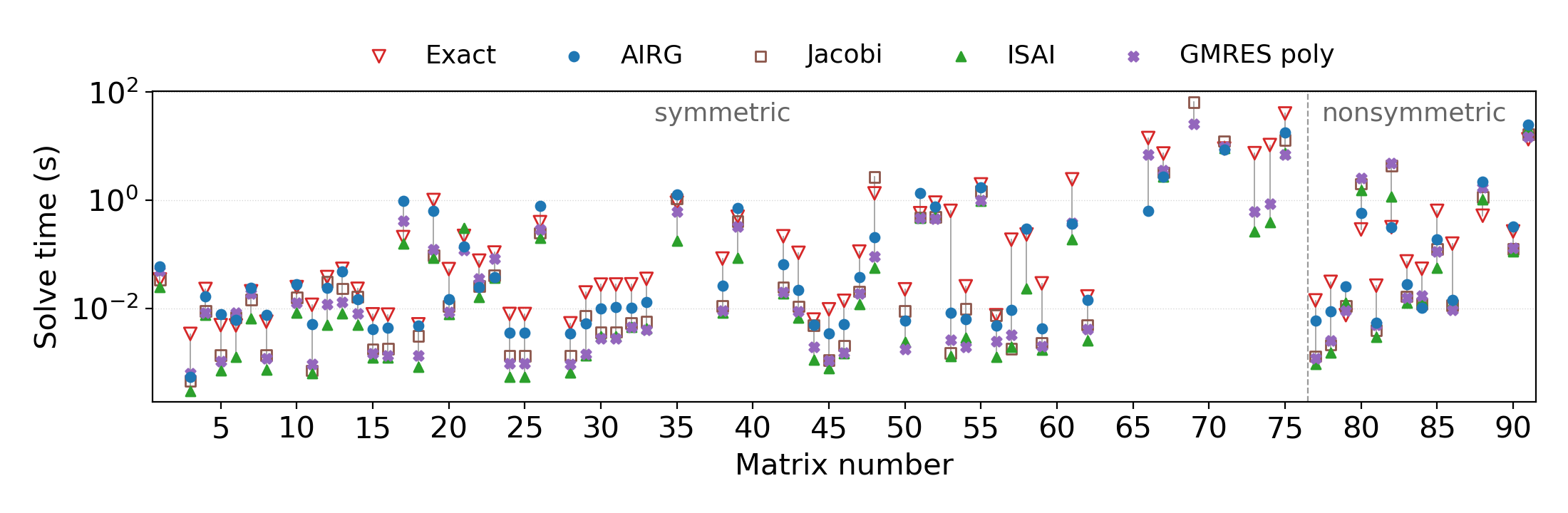}
\caption{Solve time comparison (smaller is better) for the $\mat{A}\mat{x}=\mat{b}$ solve with the ILU(1) factorisation for each test matrix, ordered as in the results tables (symmetric and nonsymmetric matrices separated by the dashed line). The solve uses an outer GMRES(30) preconditioned by the ILU factorisation, with the inner application of $\mat{L}^{-1}$ and $\mat{U}^{-1}$ performed by an exact triangular solve, by one V-cycle of AIRG, or by the Jacobi method, ISAI or the GMRES polynomial, each with the number of inner iterations that gave the smallest converged solve time for that problem. Only solves that converged to a relative tolerance of $10^{-6}$ are shown. \Cref{tab:axb_summary} of the article summarises these results.}
\label{fig:axb_timing_ilu1}
\end{figure}

\Cref{fig:setup_ilu1} is the ILU(1) counterpart of \cref{fig:setup} of the article (which shows the ILU(0) factorisation), discussed in \cref{sec:lu_solves} of the article. \Cref{fig:axb_timing_ilu0,fig:axb_timing_ilu1} give the $\mat{A}\mat{x}=\mat{b}$ solve time for each test matrix, summarised in \cref{tab:axb_summary} and discussed in \cref{sec:axb_solves} of the article.

\end{document}